\documentclass[11pt, reqno]{amsart}

\usepackage{mathrsfs,bbold}

\usepackage{amsmath, amsthm, amssymb}
\usepackage{enumerate}
\usepackage{dsfont}
\usepackage{color}
\usepackage[a4paper]{geometry}
\usepackage[utf8]{inputenc}
\usepackage{stmaryrd}
\usepackage{titlesec, fancyhdr}
\usepackage{mathabx}
\usepackage{appendix}
\usepackage{todonotes}
\usepackage{hyperref}

   \newcommand{\VV}{\mathbb{V}}

 \newcommand{\mcC}{\mathcal{C}}  \newcommand{\mcF}{\mathcal{F}}

     \newcommand{\mcT}{\mathcal{T}}

    \newcommand{\bfX}{{\bf X}}      
\newcommand{\bfV}{{\bf V}}

\renewcommand{\leq}{\leqslant}
\renewcommand{\geq}{\geqslant}

\newcommand{\ssk}{\smallskip}

\newcommand{\E}{\mathbb{E}}
\renewcommand{\P}{\mathbb{P}}

\newcommand{\R}{\mathbb{R}}
\newcommand{\N}{\mathbb{N}}
\newcommand{\Z}{\mathbb{Z}}
\newcommand{\Q}{\mathbb{Q}}

\newcommand{\vertiii}[1]{{\left\vert\kern-0.25ex\left\vert\kern-0.25ex\left\vert #1 
    \right\vert\kern-0.25ex\right\vert\kern-0.25ex\right\vert}}

\newenvironment{Dem}[1][\unskip]{
    \begin{list}
    {\hspace{0.5cm}{\sf \textbf{Proof #1 --}}}
    {   \setlength{\topsep}{0pt}%
        \setlength{\leftmargin}{0pt}%
        \setlength{\rightmargin}{0pt}%
        \setlength{\listparindent}{0pt}%
        \setlength{\itemindent}{0pt}%
        \setlength{\parsep}{0pt}%
        \addtolength{\leftmargin}{20pt}%
        \addtolength{\rightmargin}{0pt}%
	} 
	\item }{\hfill $\rhd$
	\end{list}}

\newtheoremstyle{mystyle}
{3pt}               
{3pt}               
{\it }                      
{}                      
{\sffamily\bfseries}             
{}                      
{0.5em}                 
{\llap{#2. }#1{$\;$ --}}

\theoremstyle{mystyle}

\newtheorem{thm}{Theorem}
\newtheorem*{thm*}{Theorem}

\newtheorem{cor}[thm]{\hspace{-0.15cm}  {Corollary} }
\newtheorem{lem}[thm]{\hspace{-0.2cm}  {Lemma} }
\newtheorem{prop}[thm]{\hspace{-0.2cm} {Proposition}}
\newtheorem{defn}[thm]{ \hspace{-0.3cm} {Definition}}
\newtheorem*{defn*} {Definition}
\newtheorem*{prop*} {Proposition}
\newtheorem*{lem*} {Lemma}
\newtheorem*{cor*} {Corollary}

\newtheorem{example}[thm]{\hspace{-0.15cm} {Example}}

\newtheoremstyle{mystyle2}
{3pt}               
{3pt}               
{\it }                      
{}                      
{\sffamily\bfseries}             
{}                      
{0.5em}                 
{\llap{#2 }#1{$\;$ --}}
\theoremstyle{mystyle2}

\newtheorem*{definition*}{Definition}

\titleformat{\section}[block]
{\filcenter\normalfont\sffamily\bfseries\Large}{{\hspace{-0.7cm}}\thesection \vspace{0.3cm}}{0.75em}{}

\titleformat{\subsection}[block]
{\normalfont\sffamily\bfseries\large}{\hspace{-1cm}\thesubsection \vspace{0.3cm}}{.75em}{}
\titlespacing{\subsection}{-0pc}{1.5ex plus .1ex minus .2ex}{0pc}

\titleformat{\subsubsection}[block]
{\normalfont\sffamily\bfseries}{\hspace{-1.2cm}\thesubsubsection  \vspace{0.3cm}}{.75em}{}
\titlespacing{\subsection}{0pc}{1.5ex plus .1ex minus .2ex}{0pc}

\numberwithin{equation}{section} 

\numberwithin{subsection}{section}
\numberwithin{subsubsection}{subsection}

\let\oldtocsection=\tocsection
\let\oldtocsubsection=\tocsubsection
\let\oldtocsubsubsection=\tocsubsubsection
\renewcommand{\tocsection}[2]{\hspace{0em}\oldtocsection{#1}{#2}}
\renewcommand{\tocsubsection}[2]{\hspace{1em}\oldtocsubsection{#1}{#2}}
\renewcommand{\tocsubsubsection}[2]{\hspace{2em}\oldtocsubsubsection{#1}{#2}}
\begin{document}

\vspace*{3ex minus 1ex}
\begin{center}
\huge\sffamily{Random dynamical systems, rough paths and rough flows}
\end{center}
\vskip 5ex minus 1ex

\begin{center}
{\sf I. BAILLEUL and S. RIEDEL and M. SCHEUTZOW }
\end{center}

\vspace{1cm}

\begin{center}
\begin{minipage}{0.8\textwidth}
\renewcommand\baselinestretch{0.7} \rmfamily {\scriptsize {\bf \sc \noindent Abstract.}  
We analyze common lifts of stochastic processes to rough paths/rough drivers-valued processes and give sufficient conditions for the cocycle property to hold for these lifts. We show that random rough differential equations driven by such lifts induce random dynamical systems. In particular, our results imply that rough differential equations driven by the lift of fractional Brownian motion in the sense of Friz-Victoir \cite{FV10-2} induce random dynamical systems.
}
\end{minipage}
\end{center}

\section*{Introduction}
\label{SectionIntroduction}

Rough paths theory can be seen as a pathwise solution theory for ordinary differential equations of the form
\begin{align}\label{eqn:controlled_ODE}
\begin{split}
&\dot{y_t} = b(y_t) + \sum_{i = 1}^d \sigma_i(y_t) \dot{x}^i_t; \qquad t \in [0,T] \\
&y_0 \in \R^m
\end{split}
\end{align}
where the driving signal $x \colon [0,T] \to \R^d$ is ``rough'', by which we mean only H\"older continuous with possible small H\"older exponent, and in particular not differentiable; basic references on the subject are \cite{Lyo98, LQ02, LCL07}. One of Lyons key insights was that equation \eqref{eqn:controlled_ODE} as it stands is actually not well-posed in general. More precisely, the solution map $x \mapsto y$, which is defined for smooth paths $x$, is not closable for any reasonable topology on the respective path spaces. Instead, Lyons understood that one has to enrich the path $x$ with some ``extra information'' $\mathbb{x}$, often called \emph{L\'evy area}, which is then called a \emph{rough path} $\mathbf{x} = (x,\mathbb{x})$. It turns out that \emph{rough differential equations}
\begin{align}\label{eqn:rough_DE}
\begin{split}
&dy_t = b(y_t)\, dt + \sigma(y_t) \, d\mathbf{x}_t; \qquad t \in [0,T] \\
&y_0 \in \R^m
\end{split}
\end{align}
can be defined and solved under appropriate regularity conditions on $b$ and $\sigma$, and the solution coincides with the classical one when the driving path happens to be smooth. Moreover, one can construct a topology on the space of rough paths for which the It\^o-Lyons solution map $\mathbf{x} \mapsto y$ is continuous; this result is often referred to as Lyons' \emph{universal limit theorem}. 

\ssk

One important application for rough paths theory is the pathwise solution to stochastic differential equations. Indeed, the most interesting stochastic processes do not possess smooth sample paths, like Brownian motion or general martingales. Rough paths theory allows to solve stochastic differential equations pathwise, not relying on It\^o calculus; see for instance the monographs \cite{FV10, FH14} which emphasize the applications of rough paths theory in the field of stochastic analysis. An immediate consequence is that ``rough'' stochastic differential equations generate a continuous \emph{stochastic flow}, i.e. a collection of random homeomorphisms $(\psi_{s,t})_{s,t \in [0,T]}$ on $\R^m$ which satisfy the flow property on a set of full measure; observe here that the rough path lift of a stochastic process is defined on a set of full measure, and the flow property is a result of the deterministic theory, see \cite{LQ98}. Note that establishing the flow property for an It\^o-stochastic differential equation can cause problems since the solution is a priori defined only outside a set of zero measure which depends on the whole equation, and in particular on the initial condition. The most general stochastic flows are induced by \emph{Kunita-type} stochastic differential equations which take the form
\begin{align}\label{eqn:kunita_sde}
d \psi = V(\psi,\, dt)
\end{align}
where $V$ is a stochastic process in the space of vector fields; see e.g. Kunita's monograph \cite{Kun90}. In \cite{BR15}, two of us showed that also this type of equations can be solved pathwise when enriching the vector field $V$ with a suitable second order object, and that one can extend the theory far beyond the semimartingale setting; see also \cite{der10} and \cite{DD12} where stochastic flow solutions of equations of the form \eqref{eqn:kunita_sde} where analyzed using classical rough paths theory in infinite dimensional spaces. 

\ssk

If the probability space is a metric dynamical system and the driving process satisfies the cocycle, or helix, property (such as Brownian motion), one can ask whether the solution flow to a stochastic differential equation induces a \emph{random dynamical system}, in the sense of L.~Arnold \cite{Arn98}. The theory of random dynamical systems is one of the main tools for studying the long time behaviour of stochastic flows. For instance, Oseledets' \emph{Multiplicative Ergodic Theorem} applies and yields the existence of a \emph{Lyapunov spectrum} \cite{Arn98} which describes the spatial infinitesimal behaviour of the flow and yields locally stable or unstable manifolds. Other interesting objects to study are \emph{random attractors}, that is random sets which attract points or even subsets in the state space moving according to the flow; see e.g. \cite{CF94, CDF97, C01, Cra02, S02} and \cite{DS11} in the context of stochastic differential equations. Random attractors are formulated using the language of random dynamical systems, and their properties can be analyzed by a study of the underlying random dynamical system and, in particular, the Lyapunov spectrum; see \cite{Deb98} where the Hausdorff dimension of random attractors are studied or the recent works \cite{FGS16a, FGS16b, CGS16, Vor16} in which the question of \emph{synchronization} of a random attractor is adressed. Proving the cocycle property for an It\^o-type stochastic differential equation causes similar problems as it does for the flow property since It\^o-stochastic differential equations notoriously generate nullsets. To overcome this difficulty, one of us established so-called \emph{perfection techniques} which could then be succesfully applied in the context of classical stochastic differential equations \cite{AS95, Sch96}. 

\medskip

The present work makes a first attempt to systematically connect the two fields of \emph{rough paths theory} and \emph{random dynamical systems}. In particular, we give sufficient conditions under which a random rough differential equation
\begin{align}\label{eqn:random_rough_DE}
dY_t &= b(Y_t)\, dt + \sigma(Y_t) \, d\mathbf{X}_t(\omega); \qquad t \in [0,\infty)
\end{align}
induces a random dynamical system. Our main results imply the well-known fact that a stochastic differential equation driven by a Brownian motion induces a random dynamical system. Note that our results do \emph{not} rely on perfection techniques and are rather a consequence of the pathwise calculus. Using rough paths theory, we can also go beyond the semimartingale framework and consider more general driving processes. For instance, a two-sided \emph{fractional Brownian motion} (fBm) $B^H \colon \R \to \R^d$ with Hurst index $H \in (0,1)$ is known to satisfy the helix property. Furthermore, for $H \in (1/4,1)$, the process has a natural lift to a rough paths valued process, see e.g. \cite{FV10-2}, \cite[Chapter 15]{FV10} or \cite[Chapter 10]{FH14}. Our results imply that equation \eqref{eqn:random_rough_DE} driven by this lift of fBm induces a  random dynamical system. Using the theory of \emph{rough drivers} established in \cite{BR15}, we can even take infinitely many independent fractional Brownian motions $(B^{H;i})_{i \in \N}$ as a driver, solve the (formal) equation
\begin{align*}
dY_t &= b(Y_t)\, dt + \sum_{i = 1}^{\infty} \sigma_i(Y_t) \,\circ dB^{H;i}_t(\omega); \qquad t \in [0,\infty)
\end{align*}
and our results show that also this equation induces a random dynamical system. To summarize, the results in this work yield the existence of completely new random dynamical systems and pave the way for a systematic study of the long time behaviour of flows induced by random rough differential equations driven by very general driving signals.

\medskip

Let us now comment on the existing literature. A first attempt to study the long time behaviour of stochastic differential equations driven by fractional Brownian motion was made by Hairer and coauthors, cf. \cite{Hai05, HO07, HP11, HP13} and the following works \cite{FP14, DPT16}. The approach in these works is very different from ours, as the authors do not use the theory of random dynamical systems, and define instead another object called \emph{stochastic dynamical system} and give a description of an invariant measure within this framework. Without going too much into detail here, let us mention that this theory is closer to the usual Markovian semigroup approach to invariant measures than the theory of random dynamical systems. In particular, there is currently no analogue of the \emph{Multiplicative Ergodic Theorem} in this theory, and a study of attractors seems out of reach at this stage. In a series of works, Garrido-Atienza, Lu and Schmalfuss use a pathwise calculus to generate random dynamical systems from a class of stochastic partial differential equations driven by a finite-dimensional fractional Brownian motion \cite{GLS10, GLS15, GLS16}. The strategy is similar to ours, although they do not use rough paths theory but the Riemann-Liouville representation of fBm and fractional calculus. 

\medskip

The paper is organized as follows. In Section 1 we define rough paths and rough drivers indexed by the real numbers and collect some basic properties of the respective spaces. Section 2 discusses rough paths and drivers which satisfy the cocycle property; we will call them \emph{rough cocycles}. We give sufficient conditions for this property to hold and show in particular that the lift of Brownian motion and fractional Brownian motion obtain this property. In Section 3, we prove that rough differential equations driven by rough cocycles induce random dynamical systems. We collect some basic definitions and properties of rough path spaces and rough drivers in several Appendices for the readers convenience.

\bigskip

\subsection*{Notation and basic definitions}

\ssk

\noindent \textcolor{gray}{$\bullet$} We shall always denote by $\big(E,| \cdot |\big)$ a real Banach space, not necessarily separable. Given another Banach space $F$ and a positive regularity exponent $\gamma$, we denote by $\mathcal{C}_b^{\gamma}(F,F)$ the space of all continuous, bounded functions from $F$ to itself which are $\lfloor \gamma \rfloor$-times Fr\'echet differentiable, with bounded derivatives, and for which the $\lfloor \gamma \rfloor$-th derivative is $\big(\gamma - \lfloor \gamma \rfloor\big)$-H\"older continuous. We denote by $L\big(E,\mathcal{C}_b^{\gamma}(F,F)\big)$ the space of bounded linear functions from $E$ to $\mathcal{C}_b^{\gamma}(F,F)$. 

\ssk

One can always define the algebraic tensor product $E \otimes_a E$, but a norm on this vector space is not canonically given \cite{Rya02}. In this article, we will choose a norm on $E \otimes_a E$ which is \emph{compatible}, i.e. we assume that
\begin{align}\label{eqn:compatible_tensor_norm}
 |v \otimes w| = |w \otimes v| \qquad \text{and} \qquad |v \otimes w| \leq |v| |w| \qquad \text{for all } v, w \in E.
\end{align}
The Banach space $E \otimes E$ is defined as the closure with respect to this norm. Note that compatibility is a very natural assumption and holds, for instance, for both the \emph{projective} and the \emph{injective} tensor norm, cf. \cite[Chapter 2 and 3]{Rya02}. The same assumptions are made for higher order tensor products. For $N \geq 1$, set 
\begin{align*}
  T^{N}(E) := \R \oplus E \oplus \ldots \oplus E^{\otimes N},
\end{align*}
and
\begin{align*}
  T_1^{N}(E) := 1 \oplus E \oplus \ldots \oplus E^{\otimes N}.
\end{align*}
The maps 
$$
\pi_k \colon T^{N}(E) \to E^{\otimes k}
$$ 
are the usual projection maps defined for $k = 0,\dots,N$. For $g  \in T^{N}(E)$ set 
\begin{align*}
 | g | := \max_{k = 0,\ldots,N} \big| \pi_k(g) \big|.
\end{align*}
If $g = \big(1,g^1,\ldots,g^N\big)$ and $h = \big(1,h^1,\ldots,h^N\big)$ are elements in $T_1^N(E)$, we define $g \otimes h \in T_1^N(E)$ by setting 
\begin{align*}
   \pi_k(g \otimes h) = \sum_{i = 0}^k g^{k - i} \otimes h^i
\end{align*}
for $k = 0,\ldots, N$. It can be shown that $(T_1^N(E), \otimes)$ is a topological group with unit element $\mathbf{1} = (1,0,\ldots,0)$ -- see for instance \cite[Section 2.2.1]{LCL07}. 

\medskip

\noindent \textcolor{gray}{$\bullet$} If $I$ is a subset of $\R$ and $x \colon I \to E$ is some path in $E$, we use the notation $x_{s,t} := x_t - x_s$ for $s,t  \in I$. For $p \in [1,\infty)$, we define the \emph{$p$-variation (semi-)norm} by setting 
\begin{align*}
 \| x \|_{p-\text{var};I} := \sup_{\mathbb{D}} \left( \sum_{t_i \in \mathbb{D}} |x_{t_{i+1}} - x_{t_i}|^p \right)^{\frac{1}{p}}
\end{align*}
where the supremum ranges over all finite subsets $\mathbb{D} = \{t_0 < t_1 < \ldots < t_N \}$ of $I$. Note that this seminorm becomes a norm if we restrict it to a space of functions with value $0$ at some given time point. If $x, y \colon I \to E$ are two paths, we set 
\begin{align*}
  d_{p-\text{var};I}(x,y) := \|x - y\|_{p-\text{var};I}.
\end{align*}

\medskip

\noindent \textcolor{gray}{$\bullet$} Let $(\Omega,\mathcal{F})$ and $(X,\mathcal{B})$ be measurable spaces.
A family $\theta=(\theta_t)_{t \in \mathbb{R}}$ of maps from $\Omega$ to itself is called a \emph{measurable dynamical system} if
\begin{itemize}
   \item[(i)] $(\omega,t) \mapsto \theta_t\omega$ is $\mcF\otimes \mathcal{B}(\R) / \mcF$-measurable,   \vspace{0.07cm}

   \item[(ii)] $\theta_0 = \operatorname{Id}$,   \vspace{0.07cm}
  
   \item[(iii)] $\theta_{s + t} = \theta_s \circ \theta_t$, for all $s,t \in \mathbb{R}$.   \vspace{0.1cm}
\end{itemize}
If $\P$ is furthermore a probability on $(\Omega,\mathcal{F})$ that is invariant by any of the elements of $\theta$,
$$
\P \circ \theta_t^{-1} = \P
$$ 
for every $t \in \mathbb{R}$, we call the tuple $\big(\Omega, \mathcal{F},\P,\theta\big)$ a \emph{measurable metric dynamical system}. Let $\mathbb{T}$ be either $\R$ or $[0,\infty)$, equipped with its Borel $\sigma$-algebra. A \emph{measurable random dynamical system} on $(X,\mathcal{B})$ is a measurable metric dynamical system $\big(\Omega, \mathcal{F},\P,\theta\big)$ with a measurable map 
   $$
   \varphi \colon \mathbb{T} \times \Omega \times X \to X
   $$ 
   that enjoys the \emph{cocycle property}, i.e. $\varphi_0(\omega) = \operatorname{Id}_X$, for all $\omega \in \Omega$, and
  \begin{align}\label{eqn:cocycle_property}
   \varphi_{t+s}(\omega) = \varphi_t(\theta_s\omega) \circ \varphi_s(\omega)
  \end{align}
  for all $s,t \in \mathbb{T}$ and $\omega \in \Omega$. If $X$ is a topological space and the map $\varphi_{\cdot}(\omega,\cdot) \colon \mathbb{T} \times X \to X$ is continuous for every $\omega \in \Omega$, it is called a \emph{continuous random dynamical system}.

\bigskip

\section{Rough paths and rough drivers on the real line}
\label{SectionRealLine}

One of the main sources of continuous time random dynamical systems is given by the random flows generated by stochastic differential equations of It\^o-Stratonovich or Kunita-type; see Arnold's book \cite{Arn98} for a well-documented reference on the subject. We show in the next two section that the setting of rough flows introduced in \cite{BR15} offers a wealth of other examples beyond the classical paradigm of It\^o's stochastic calculus. To make a long story short, while T. Lyons' theory of rough paths provides a deterministic and pathwise view on stochastic differential equations, the setting of rough flows provides a deterministic counterpart of Kunita-Le Jan-Watanabe theory of stochastic flows \cite{LeJ82, Kun90}. We shall not give here an introduction to these theories and refer instead a newcomer in these fields to the pedagogic lecture notes \cite{FH14, BaiLN}, and the articles \cite{BR15, BC16}. Appendix \ref{AppendixRP} contains however some elementary definitions and facts on rough paths in a Banach setting, and Appendix \ref{AppendixRoughFlows} gives a bird's eye view on rough drivers and rough flows. Nothing more than a working knowledge of the definition of a rough path and a rough driver, and the basic well-posedness results for rough differential equations on paths and flows, will be used in the present work. We invite the reader to have in mind the image of a rough path as an integrator, playing the role of a multi-dimensional control $h$ in a controlled ordinary differential equation $\dot x_t = \sum V_i(x_t)\,\dot h^i_t$, for some given vector fields $V_i$ on $\R^m$ say, and to have in mind the image of a rough driver as a time dependent vector field. The fields $V_i(\cdot)\,\dot h^i_t$ that appear in controlled ordinary differential equations, and their generalisation with rough paths as controls, are special cases of such rough drivers.

\medskip

Rough paths and rough drivers are usually defined on compact intervals; the extension of these notions on  the whole real line is done as follows. Set
\begin{align*}
\mathcal{I}_0 := \Big\{ I \subset \R\, :\, I \text{ is a compact interval containing } 0 \Big\}.
\end{align*}

\medskip

\begin{definition*} {\normalfont
	We denote by $\mathcal{C}_0^{0,p-\text{var}}\big(\R,T_1^N(E)\big)$ the space of continuous paths 
	$$
	\mathbf{x} \colon \R \to T_1^N(E)
	$$ 
	for which $\mathbf{x} \vert_I \in \mathcal{C}_0^{0,p-\text{var}}\big(I,T_1^N(E)\big)$, for every interval $I \in \mathcal{I}_0$. We equip this set with the coarsest topology for which the projection maps
	\begin{align*}
	\mathcal{C}_0^{0,p-\text{var}}\big(\R,T_1^N(E)\big) &\to \mathcal{C}_0^{0,p-\text{var}}\big(I,T_1^N(E)\big)   \\
	\mathbf{x} &\mapsto \mathbf{x} \vert_{I}
	\end{align*}
	are continuous for all  intervals $I \in \mathcal{I}_0$.
	}
\end{definition*}

\medskip

\begin{prop}
Let $p \in [1,\infty)$ and $N \geq 1$. The space $\mathcal{C}_0^{0,p-\text{var}}\big(\R,T_1^N(E)\big)$ is a completely metrizable topological space. It is separable, and therefore Polish, if and only if $E$ is separable.
\end{prop}

\medskip  

\begin{Dem}
One can easlily check that
	\begin{align*}
	d_{p-\text{var}}(\mathbf{x},\mathbf{y}) := \sum_{m = 1}^{\infty} 2^{-m} \Big( d_{p-\text{var};[-m,m]}(\mathbf{x},\mathbf{y}) \wedge 1 \Big)
	\end{align*}
is a complete metric which induces the topology on the space. It remains to consider separability. We can extend every continuous path $\mathbf{x} \colon [0,1] \to T_1^N(E)$ to a continuous path on the whole real line by setting it constant outside $[0,1]$. Thus, we see that the space $\mathcal{C}_0^{0,p-\text{var}}\big([0,1],T_1^N(E)\big)$ is homeomorphic to a subset of $\mathcal{C}_0^{0,p-\text{var}}\big(\R,T_1^N(E)\big)$. If $E$ is not separable, we know from Proposition \ref{prop:rp_spaces_comp_Polish} that $\mathcal{C}_0^{0,p-\text{var}}\big([0,1],T_1^N(E)\big)$ is not separable which implies that $\mathcal{C}_0^{0,p-\text{var}}\big(\R,T_1^N(E)\big)$ can not be separable either. Next, assume that $E$ is separable. 	From Proposition \ref{prop:rp_spaces_comp_Polish}, we know that $\mathcal{C}_0^{0,p-\text{var}}\big(I,T_1^N(E)\big)$ is a Polish spaces for every interval $I \in \mathcal{I}_0$. Set $X_n := \mathcal{C}_0^{0,p-\text{var}}\big([-n,n],T_1^N(E)\big)$ for $n \in \N$. If $n \geq m$, we define maps 
$\varphi_{n m} \colon X_n \to X_m$ by setting $\varphi(\mathbf{x}) = \mathbf{x}\vert_{[-m,m]}$. One can check that $\{X_n,\varphi_{n m},\N \}$ defines a projective system of topological spaces and that the projective limit is homeomorphic to $\mathcal{C}_0^{0,p-\text{var}}\big(\R,T_1^N(E)\big)$. Recall that the projective limit is defined as a subset of the product space $\Pi_{n = 1}^{\infty} X_n$, and that this subset is even a closed subspace of the product space since all $X_n$ are Hausdorff spaces -- see for instance \cite[Proposition 1.1.1 and Lemma 1.1.2]{RZ10}. Since countable products and closed subspaces of Polish spaces are again Polish \cite[Proposition 8.1.3. and Proposition 8.1.1.]{Coh80}, the projective limit is a Polish space and the claim follows.
\end{Dem}

\medskip

We have a similar definition for rough drivers.

\medskip

\begin{definition*} {\normalfont
Let $2\leq p<3$ and $p-2<\rho\leq 1$ be given. A weak (geometric) $(p,\rho)$-rough driver on the real line $\R$ is a family of maps ${\bfV}_{s,t} = \big(V_{s,t},\VV_{s,t}\big)$, $s \leq t \in \R$, such that all ${\bfV} \vert_{I} = \big(V\vert_{I},\VV \vert_{I}\big)$ are (geometric) $(p,\rho)$-rough driver for every $I \in \mathcal{I}_0$. If $\bfV$ and $\bfV'$ are (geometric) $(p,\rho)$-rough drivers, set 
 \begin{align*}
  d_{p,\rho}(\bfV,\bfV') := \sum_{m = 1}^{\infty} 2^{-m} \Big( d_{p,\rho;[-m,m]}(\bfV \vert_{[-m,m]}, \bfV' \vert_{[-m,m]}) \wedge 1 \Big).
 \end{align*}
We denote by $\mathcal{D}^{p,\rho}_g$ the set of geometric $(p,\rho)$-rough drivers equipped with the metric $d_{p,\rho}$.   }
\end{definition*}

\bigskip

\section{Rough cocycles}
\label{SectionRoughCocycle}

Working in the setting of rough paths and rough drivers, our building blocks for constructing continuous time random dynamical systems will thus be It\^o-type and Kunita-type cocycles, such as defined here. 

\medskip

\begin{defn}\label{def:p_rough_helix} {\normalfont
Let $\big(\Omega, \mathcal{F},\P,\theta\big)$ be a measurable metric dynamical system.   \vspace{0.15cm}

\begin{itemize}
   \item[\textcolor{gray}{$\bullet$}] Let $p \in [1,\infty)$. A process $\mathbf{X} \colon \R \times \Omega \to T_1^{\lfloor p \rfloor}(E)$ is called a geometric \emph{$p$-rough path cocycle} if  
$$
\mathbf{X}(\omega) \in \mathcal{C}_0^{0,p-\text{var}}\big(\R,T_1^{\lfloor p \rfloor}(E)\big),\quad \textrm{ for every } \omega \in \Omega,
$$ 
and if the cocycle relation
\begin{align}\label{eqn:helix_property_1}
  \mathbf{X}_{s+t}(\omega) = \mathbf{X}_{s}(\omega) \otimes \mathbf{X}_{t}(\theta_s\omega)
\end{align}
holds for every $\omega \in \Omega$ and every $s,t \in \R$, that is 
\begin{align}\label{eqn:helix_property_2}
\mathbf{X}_{s, s + t}(\omega) = \mathbf{X}_{t}(\theta_s\omega),
\end{align}
with the increment notation.   \vspace{0.1cm}

   \item[\textcolor{gray}{$\bullet$}] Let $2\leq p<3$ and $p-2<\rho\leq 1$ be given. A $\mathcal{D}^{p.\rho}_g$- valued random variable $\bfV$ is called a geometric $p$-\emph{rough driver cocycle} if one has the identity
\begin{equation}
\label{EqCocycleRoughDriver}
   \bfV_{s,s+t}(\omega) = \bfV_{0,t}(\theta_s\omega)   
\end{equation}
for every $\omega \in \Omega$ and every $s \in \R$ and $t \in [0,\infty)$.
\end{itemize}
  }
\end{defn}

\medskip

We simply talk of a \textit{rough path} or a \textit{rough driver cocycle}.

\bigskip

\subsection{Existence of rough path cocycles}
\label{Subsectionrpcocycles}

Our aim in this section is to give sufficient conditions under which a rough paths valued process defines a rough path cocycle. We start with a simple identity for iterated integrals of smooth paths. For a path $x \colon \R \to E$, we set as usual $x_{s,t} := x_t - x_s$, and define
\begin{align*}
(\vartheta_u x)_t := x_{t+u} - x_u. 
\end{align*}

\medskip

\begin{lem}\label{lem:riem_sums_calc}
	Let $x, y \colon \R \to E$ be two $E$-valued paths. Let $[u,v]$ be a closed interval in $\R$ and $u = t_0 < \ldots < t_N = v$ be a finite partition of $[u,v]$. Choose $t'_i \in [t_i,t_{i+1}]$ for $i = 0,\ldots,N-1$. Then for every $h \in \R$,
	\begin{align*}
	\sum_{i = 0}^{N-1} x_{s,t'_i} \otimes y_{t_i,t_{t+1}} = \sum_{i = 0}^{N-1} (\vartheta_h x)_{s-h,\hat{t}'_{i}} \otimes (\vartheta_h y)_{\hat{t}_i, \hat{t}_{i+1}}
	\end{align*}
	where $\hat{t}_i = t_i - h$ and $\hat{t}'_i = t'_i - h$. In particular, provided the Riemann sums above converge, we have the identity
	\begin{align*}
	\int_s^t x_{s,u} \, \otimes dy_u = \int_{s-h}^{t-h} (\vartheta_h x)_{s-h,u} \, \otimes d(\vartheta_h y)_u.
	\end{align*}
\end{lem}

\medskip

\begin{Dem}
	This is just a simple calculation, since
	\begin{align*}
	&\sum_{i} \big(x(t'_i) - x(s)\big) \otimes \big(y(t_{i+1}) - y(t_i)\big)   \\
	&= \sum_{i} \big((\vartheta_h x)(t'_i - h) + x(h) - x(s)\big) \otimes \big((\vartheta_h y)(t_{i+1} - h) - (\vartheta_h y)(t_{i} - h)\big)   \\
	&= \sum_{i} \big((\vartheta_h x)(\hat{t}'_i ) - (\vartheta_h x)(s - h)\big) \otimes \big((\vartheta_h y)(\hat{t}_{i+1} ) - (\vartheta_h y)(\hat{t_{i}} )\big).
	\end{align*}
\end{Dem}

\medskip

\begin{prop}
	\label{prop:smooth_rp_are_cocycles}
	Let $x \colon \R \to E$ be a continuous path which is locally of finite $p$-variation for some $p \in [1,2)$. Then
	\begin{align*}
	\int_{\Delta^n_{s,t}} dx \otimes \cdots \otimes dx = \int_{\Delta^n_{s-h,t-h}} d(\vartheta_h x) \otimes \cdots \otimes d(\vartheta_h x) \in E^{\otimes n}
	\end{align*}
	for every interval $[s,t]$ and every $h \in \R$ where the integrals exist as Young integrals.
\end{prop}

\medskip

\begin{Dem}
	The Young integrals exist as limit of Riemann sums -- see for instance \cite[Theorem 1.16]{LCL07} for a proof for Banach space valued paths. The equality follows by applying Lemma \ref{lem:riem_sums_calc} repeatedly.
\end{Dem}

\bigskip

Although elementary, the next result will be our workhorse in this section.

\medskip

\begin{thm}\label{prop:can_rough_cocycle}
Let $p \in [1,\infty)$ and $N \geq 1$ be given. Let also $\big(\,\overline{\Omega}, \overline{\mathcal{F}},\overline{\P}\,\big)$ be a probability space and let $\overline{\mathbf{X}}$ be a $\mathcal{C}_0^{0,p-\text{var}}\big(\R,T_1^N(E)\big)$-valued random variable defined on that probability space. Assume that $\overline{\mathbf{X}}$ has stationary increments, that is the distribution of the process $\left( \overline{\mathbf{X}}_{t_0,t_0+h}\right)_{h \in \R}$ does not depend on $t_0 \in \R$. Then there exists a metric dynamical system $\big(\Omega,\mathcal{F},\P, \theta\big)$ and a $\mathcal{C}_0^{0,p-\text{var}}\big(\R,T_1^N(E)\big)$-valued random variable ${\mathbf{X}}$ defined on $\Omega$, such that $\mathbf{X}$ has the same law as $\overline{\mathbf{X}}$, while $\mathbf{X}$ enjoys the cocycle property \eqref{eqn:helix_property_1}, or \eqref{eqn:helix_property_2}.
\end{thm}

\medskip

\begin{Dem}
Set
\begin{align*}
\Omega := \mathcal{C}_0^{0,p-\text{var}}\big(\R,T^N_1(E)\big),
\end{align*}
and let $\mathcal{F}$ be the Borel $\sigma$-algebra and $\P$ be the law of $\overline{\mathbf{X}}$. Define $\theta$ by
\begin{align*}
(\theta_t \omega)(s) = \omega(t)^{-1} \otimes \omega(t+s)
\end{align*}
and define $\bfX$ as the canonical coordinate process 
$$
\mathbf{X}_t(\omega) := \omega(t).
$$ 
Since $\overline{\mathbf{X}}$ has stationary increments, $\P$ is invariant under any shift $\theta_t$, for every $t \in \R$. The only property we need to prove is the joint measurability of the map $(t, \omega) \mapsto \theta_t \omega$. Recall the definition of the rough path-lift operator on 'smooth' path, given in  equation \eqref{EqDefnLift}, in Appendix \ref{AppendixRP}. By definition, we know that there is a dense subset $\Omega_0 \subset \Omega$ such that every $\omega \in \Omega_0$ has the property that for every interval $(-a,b)$ containing $0$,
\begin{align*}
\omega \vert_{(-a,b)} = S_N(x)
\end{align*}
for some $x \in \mathcal{C}_0^{0,1-\text{var}}\big((-a,b),E\big)$. In a first step, we show that $t \mapsto \theta_t \omega$ is continuous for every $\omega \in \Omega_0$. Fix $\omega \in \Omega_0$ and let $t_n \to t$. We have to prove that $d_{p-\text{var};(-a,b)}\big(\theta_{t_n} \omega, \theta_t \omega\big) \to 0$ for $n \to \infty$. Choose $m$ sufficiently large such that $(s + t_n) \in [-m,m]$ for every $s \in (-a,b)$ and every $n \in \N$. Let $x \in \mathcal{C}_0^{0,1-\text{var}}\big([-m,m],E\big)$ be such that $S_N(x) = \omega \vert_{[-m,m]}$. By continuity of $S_N$, the claim follows from Proposition \ref{prop:smooth_rp_are_cocycles} if we can prove that
\begin{align*}
\big\| \vartheta_t x - \vartheta_{t_n} x \big\|_{1-\text{var};I(-a,b)} \to 0
\end{align*}
for $n \to \infty$ where $\vartheta_t x = x_{t + \cdot} - x_t$ . By definition of the space $\mathcal{C}_0^{0,1-\text{var}}\big([-m,m],E\big)$ as the closure of smooth paths, we may assume that $x$ itself is smooth. In this case,
\begin{align*}
\big\| \vartheta_t x - \vartheta_{t_n} x \big\|_{1-\text{var};(-a,b)} = \int_{-a}^b \big|x'_{t+s} - x'_{t_n + s}\big|\, ds \to 0
\end{align*}
for $n \to \infty$ by continuity of $x'$ and Lebesgue's dominated convergence theorem; so $t \mapsto \theta_t \omega$ is indeed continuous for $\omega \in \Omega_0$. 

\ssk

Now let $\omega \in \Omega$ be arbitrary and $\varepsilon > 0$ be given. Since $\Omega_0$ is dense, we can find an $\hat{\omega} \in \Omega_0$ such that $d_{p-\text{var};[-m,m]}(\omega, \hat{\omega}) \leq \varepsilon/3$. By the triangle inequality, $d_{p-\text{var};(-a,b)}\big(\theta_{t_n} \omega, \theta_t \omega\big)$ is no greater than 
\begin{align*}
&d_{p-\text{var};(-a,b)}\big(\theta_{t_n} \omega, \theta_{t_n} \hat{\omega}\big) + d_{p-\text{var};(-a,b)}\big(\theta_{t_n} \hat{\omega}, \theta_t \hat{\omega}\big) + d_{p-\text{var};(-a,b)}\big(\theta_{t} \hat{\omega}, \theta_t \omega\big)   \\
&\leq 2 d_{p-\text{var};[-m,m]}(\omega, \hat{\omega}) + d_{p-\text{var};(-a,b)}\big(\theta_{t_n} \hat{\omega}, \theta_t \hat{\omega}\big) \leq \varepsilon
\end{align*}
for $n$ sufficiently large. This shows that $t \mapsto \theta_t \omega$ is continuous for every fixed $\omega \in \Omega$. Now fix $t \in \R$ and assume $\omega_n \to \omega$. Choose $m \geq 1$ such that $s + t \in [-m,m]$ for every $s \in (-a,b)$; then
\begin{align*}
d_{p-\text{var};(-a,b)}\big(\theta_{t} \omega_n, \theta_t \omega\big) \leq d_{p-\text{var};[-m,m]}\big(\omega_n, \omega\big) \to 0
\end{align*}
for $n \to \infty$, thus $\omega \mapsto \theta_t \omega$ is continuous for every fixed $t \in \R$ and in particular measurable. This implies that $(t, \omega) \mapsto \theta_t \omega$ is jointly measurable.
\end{Dem}

\medskip

\begin{cor}\label{cor:coc_prop_stable_weak_conv}
Let $(\mathbf{X}^n)_{n\geq 1}$ be a sequence of stochastic processes with values in $\mathcal{C}_0^{0,p-\text{var}}\big(\R,T_1^N(E)\big)$, and stationary increments. Assume they converge in law to some limit process $\overline{\mathbf{X}}$. Then there exists a metric dynamical system $\big(\Omega,\mathcal{F},\P, \theta\big)$ and a $\mathcal{C}_0^{0,p-\text{var}}\big(\R,T_1^N(E)\big)$-valued random variable ${\mathbf{X}}$ with the same law as $\overline{\mathbf{X}}$, which enjoys furthermore the cocycle property \eqref{eqn:helix_property_2}.
\end{cor}

\medskip

\begin{Dem}
 Using that the laws of the processes $\mathbf{X}^n$ and $\overline{\mathbf{X}}$ are determined through their finite dimensional distributions, one sees that $\overline{\mathbf{X}}$ also has stationary increments, so the claim follows from Theorem \ref{prop:can_rough_cocycle}.
\end{Dem}

\medskip

In the following, we show how one can use natural approximation procedures to prove that some interesting classes of random processes enjoy the cocycle property.

\bigskip

\subsubsection{Approximation by convolution}
\label{SubsubsectionConvolution}

Convolution provides a very natural approximation of a stochastic process which satisfies the cocycle property because the latter is preserved, as this section will make it clear. Let $X \colon \R \times \Omega \to E$ be a continuous stochastic process, i.e. $X$ is jointly measurable and $t \mapsto X_t(\omega)$ is continuous for every $\omega \in \Omega$. Let $\mu$ be a measure on $\R$ with a continuous, compactly supported density. Set
\begin{equation*}
X^{\mu}_t(\omega) := \int \Big(X_{t - u}(\omega) - X_{- u}(\omega)\Big)\, \mu(du).
\end{equation*}
Note that the integral indeed exists as a limit of Riemann sums which implies that $X^{\mu} \colon \R \times \Omega \to E$ is again a stochastic process with continuous trajectories. A straightforward calculation justifies the following fact.

\medskip

\begin{lem}\label{lem:coc_prop_stable_conv}
Let $\big(\Omega,\mathcal{F},\P,\theta\big)$ be a measurable metric dynamical system and $X \colon \R \times \Omega \to E$ a continuous stochastic process which satisfies the cocycle property. Then $X^{\mu}$ also satisfies the cocycle property.
\end{lem}

\medskip

Our main result concerning approximations by convolution is the following.

\medskip

\begin{prop}  \label{prop:rough cocylce_convolution}
Let $X \colon \R \to E$ be a continuous stochastic process with stationary increments and $X_0 = 0$ almost surely. Let $\mu_n$ be a sequence of compactly supported probability measures on $\R$, with a continuously differentiable density with respect to Lebesgue measure, and such that $\mu_n$ converges weakly to the Dirac mass $\delta_0$ as $n$ tends to $\infty$. Denote by $\mathbf{X}^{n}$ the canonical lift of $X^{\mu_n}$ to the space $\mathcal{C}_0^{0,p-\text{var}}(\R,T_1^N(E))$. Assume that the $\mathbf{X}^{n}$ converge in law to some limit process $\overline{\mathbf{X}}$. Then there exists a metric dynamical system $\big(\Omega,\mathcal{F},\P, \theta\big)$ and a $\mathcal{C}_0^{0,p-\text{var}}(\R,T_1^N(E))$-valued random variable ${\mathbf{X}}$ with the same law as $\overline{\mathbf{X}}$, while it also enjoys the cocycle property \eqref{eqn:helix_property_2}.
\end{prop}

\medskip

\begin{Dem}
We may assume, without loss of generality, that the underlying probability space $(\Omega,\mathcal{F},\P)$ is given by $\Omega = \mathcal{C}^0(\R,E)$, with $\mathcal{F}$ its Borel $\sigma$-algebra and $\P$ is the probability measure for which $X$ is given as the canonical process. Moreover, we may assume that $X$ satisfies the cocycle property for the shift map $\theta_t \omega = \omega(\cdot + t) - \omega(t)$. By Lemma \ref{lem:coc_prop_stable_conv}, it follows that every $X^{\mu_n}$ satisfies the cocycle property. Since the sample paths of the processes $X^{\mu_n}$ are smooth, we can apply Proposition \ref{prop:smooth_rp_are_cocycles} and see that the canonical lifts $\mathbf{X}^{n}$ of $X^{\mu_n}$ also satisfy the cocycle property for every $n \in \N$. By assumption, $\mathbf{X}^{n}$ converges weakly to $\overline{\mathbf{X}}$ as $n$ goes to $\infty$, so the claim follows from Corollary \ref{cor:coc_prop_stable_weak_conv}.
\end{Dem}

\medskip

It follows that the natural rough paths lift of a Gaussian process with stationary increments defines a cocycle. Recall the notion of finite $\varrho$-variation for a function indexed by a $2$-dimensional parameter, such as used in rough paths theory in the setting of Gaussian rough paths \cite{FV10-2, FH14}.

\medskip

\begin{cor}\label{cor:cocycle_gaussian}
Let $X \colon \R \to \R^d$ be a continuous, centered Gaussian process with stationary increments and independent components. Assume that its covariance function has finite $2$-dimensional $\varrho$-variation on every square $[s,t]^2$ in $\R^2$, for some $\varrho \in [1,2)$. Let $\overline{\mathbf{X}}$ stand for the natural lift of $X$, in the sense of Friz-Victoir, with sample paths in the space $\mathcal{C}_0^{0,p-\text{var}}\big(\R,T_1^{\lfloor p \rfloor}(\R^d)\big)$, for every $p > 2\varrho$. There exists a metric dynamical system $\big(\Omega,\mathcal{F},\P, \theta\big)$ and a $\mathcal{C}_0^{0,p-\text{var}}\big(\R,T_1^{\lfloor p \rfloor}(\R^d)\big)$-valued random variable ${\mathbf{X}}$ with the same law as $\overline{\bfX}$ which also enjoys the cocyle property.
\end{cor}

\medskip

\begin{Dem}
The assumption of finite $\varrho$-variation of the covariance function is the key assumption which guarantees the existence of a natural lift $\overline{\mathbf{X}}$ on every compact interval $[-L,L]$, $L \geq 1$ -- cf. \cite[Theorem 15.33]{FV10} and \cite{FGGR16} for a further discussion. Glueing together these lifts, we obtain a process $\overline{\mathbf{X}}$ which is defined on the whole real line. Note that on compact intervals, it follows from \cite[Theorem 15.45]{FV10} that $\overline{\mathbf{X}}$ is the limit in probability of approximations by convolutions $\mathbf{X}^n$ in the sense above, so we also have the convergence
\begin{align*}
\mathbf{X}^n \to \overline{\mathbf{X}}
\end{align*}
in probability for $n\to \infty$, and we can conclude with Proposition \ref{prop:rough cocylce_convolution}.
\end{Dem}

\medskip

\begin{example} {\normalfont
	The result above can be applied to the lift of a (two-sided) fractional Brownian motion $B^H  \colon \R \to \R^d$ with Hurst parameter $H > 1/4$ for $p > 1/H$ (\cite[Proposition 15.5]{FV10}).
	}
\end{example}

\bigskip

\subsubsection{Piecewise linear approximations}
\label{SubsubsectionPiecewiseApprox}

Let $X \colon \R \to E$ be a stochastic process and 
$$
\mathbb{D} = \{ \ldots < t_i < t_{i+1} < \ldots \}
$$
be a countable ordered subset of $\R$ for which $\inf \mathbb{D} = - \infty$ and $\sup \mathbb{D} = \infty$. We define the piecewise linear approximation $X^\mathbb{D}$ of $X$ with respect to $\mathbb{D}$ by setting
\begin{align*}
 X^\mathbb{D}_t := X_{t_i} + (t - t_i) \frac{X_{t_{i+1}} - X_{t_{i}}}{t_{i+1} - t_i}
\end{align*}
for $t \in [t_i,t_{i+1}]$. It should be clear that we cannot expect the cocycle property to hold for $X^{\mathbb{D}}$ even if it holds for $X$. However, the following weak form still holds.

\medskip

\begin{lem}\label{lem:coc_prop_alm_stable_piecewise_lin}
Let $\big(\Omega,\mathcal{F},\P,\theta\big)$ be a measurable metric dynamical system and $X \colon \R \times \Omega \to E$ be a stochastic process enjoying the cocycle property. Assume that $\mathbb{D}$ consists of equidistant consecutive points, so $t_{i+1} - t_i =: \delta$ does not depend on $i$. Then
\begin{align*}
X^\mathbb{D}_{t}(\theta_h \omega) - X^\mathbb{D}_{s}(\theta_h \omega) = X^\mathbb{D}_{t + h}(\omega) - X^\mathbb{D}_{s + h}(\omega)   
\end{align*}
for every $s,t \in \R$ and every $h \in \delta \Z$.
\end{lem}

\medskip

Similarly as for the convolution, we have the following result.

\medskip

\begin{prop}\label{prop:rough cocylce_dyadic}
Let $X \colon \R\times\overline\Omega \to E$ be a continuous stochastic process with stationary increments and $X_0 = 0$ almost surely. Set $\mathbb{D}_n := \{k 2^{-n}\, :\, k \in \Z \}$, and denote by $\mathbf{X}^{n}$ the canonical lift of $X^{\mathbb{D}_n}$ to the space $\mathcal{C}_0^{0,p-\text{var}}\big(\R,T_1^N(E)\big)$, where $p \geq 1$. If ${\bfX}^n$ converges in law to some limit process $\overline{\bfX}$, then there exists a metric dynamical system $\big(\Omega,\mathcal{F},\P, \theta\big)$ and a $\mathcal{C}_0^{0,p-\text{var}}\big(\R,T_1^N(E)\big)$-valued random variable $\mathbf{X}$ with the same law as $\overline{\bfX}$ that has furthermore the cocyle property \eqref{eqn:helix_property_2}.
\end{prop}

\medskip

\begin{Dem}
The proof is similar to the proof of Proposition \ref{prop:rough cocylce_convolution}. Again, we may assume that the underlying probability space $(\Omega,\mathcal{F},\P)$ is given by $\Omega = \mathcal{C}^0(\R,E)$, with $\mathcal{F}$ its Borel $\sigma$-algebra and $\P$ is the probability measure for which $X$ is given as the canonical process. Set $\mathbb{D} := \bigcup_n \mathbb{D}_n$. Using Lemma \ref{lem:riem_sums_calc} and Lemma \ref{lem:coc_prop_alm_stable_piecewise_lin}, we can conclude that
\begin{align}\label{eqn:dyadic_indep_shift}
\mathcal{L}\big(\,\overline{\mathbf{X}}_{t,t+h}\, :\, h \in \R\big)
\end{align}
does not depend on $t \in \mathbb{D}$. Since the process $\overline{\mathbf{X}}$ is continuous, the family of laws \eqref{eqn:dyadic_indep_shift} actually does not depend on $t \in \R$, and the claim follows from Proposition \ref{prop:can_rough_cocycle}.
\end{Dem}

\medskip

\begin{cor}
Let $M \colon \R \to \R^d$ be a continuous semimartingale with stationary increments. Denote by $\overline{\mathbf{M}} \colon \R \to T_1^2(\R^d)$ the Stratonovich lift of $M$ for some $p \in (2,3)$. Then there exists a metric dynamical system $\big(\Omega,\mathcal{F},\P, \theta\big)$ and a $\mathcal{C}_0^{0,p-\text{var}}\big(\R,T_1^{2}(\R^d)\big)$-valued random variable ${\mathbf{M}}$ which has the same law as $\overline{\mathbf{M}}$ and enjoys the cocycle property.
\end{cor}

\medskip

\begin{Dem}
We know from \cite[Theorem 14.16]{FV10} that $\overline{\mathbf{M}}$ can be approximated in $p$-variation metric on compact sets by the canonical lift of dyadic approximations in probability. (Note that \cite[Theorem 14.16]{FV10} is actually formulated for continuous local martingales only. However, on \cite[p. 386]{FV10} it is explained why this result implies the same statement for semimartingales.) We conclude as in Corollary \ref{cor:cocycle_gaussian} using Proposition \ref{prop:rough cocylce_dyadic}.
\end{Dem}

\medskip

In \cite{LLQ02}, Ledoux, Lyons and Qian consider rough path lifts of $V$-valued Wiener processes. We quickly recall the definition here. Assume that $\big(E,\mathcal{H},\gamma\big)$ is a \emph{Gaussian Banach space}, that is $E$ is a separable Banach space, $\gamma$ is a centered Gaussian measure on the Borel sets of $E$ and $\mathcal{H}$ denotes the \emph{Cameron-Martin space}, that is a separable Hilbert space which is continuously embedded in $E$. Note that by the Riesz representation theorem, the dual space $E^*$ can be considered as a subspace of $\mathcal{H}$. \emph{Gaussian measure} means that every continuous, linear functional $\lambda$ on $E$ has a Gaussian law $\mathcal{N}\big(0,|\lambda|_{\mathcal{H}}^2\big)$. Let $\big(\,\overline{\Omega}, \overline{\mathcal{F}}, \overline{\P}\,\big)$ be a probability space. A stochastic process $X \colon \R \to E$ is called \emph{Wiener process based on $\big(E,\mathcal{H},\gamma\big)$} if
\begin{enumerate}
 \item $X_0 = 0$ almost surely,   \vspace{0.1cm}
 
 \item $X$ has independent increments and
 \begin{align*}
  \lambda(X_t - X_s) \sim \mathcal{N}\Big(0, (t-s) |\lambda|^2_{\mathcal{H}}\Big)
 \end{align*}
 for all $s \leq t$, and $\lambda \in E^*$,   \vspace{0.1cm}
 
 \item almost all sample paths $t \mapsto X_t(\omega)$ are continuous.   \vspace{0.1cm}
\end{enumerate}

For proving the existence of a natural lift of such a process, Ledoux, Lyons and Qian introduce a further condition on the norm chosen on the algebraic tensor product $E \otimes_a E$ and the Gaussian measure $\gamma$. If this condition is satisfied, they call the pair $\big(| \cdot |_{E \otimes E}, \gamma\big)$ \emph{exact}. In order not to distract the reader with technical details, we do not to repeat the definition of exactness here and refer her/him instead to \cite[Definition 1, p. 565]{LLQ02}. If this condition is satisfied, they can show that for each fixed $(-a,b)$, the natural lifts $\mathbf{X}^n$ of the dyadic approximations $X^n$ converge in the space $\mathcal{C}_0^{0,p-\text{var}}\big((-a,b),T_1^2(E)\big)$ almost surely for every $p \in (2,3)$ -- \cite[Theorem 2 and 3]{LLQ02}. We can glue the processes together to obtain a lift $\overline{\mathbf{X}}$ which has almost surely sample paths in $\mathcal{C}_0^{0,p-\text{var}}\big(\R,T_1^2(E)\big)$. Using Proposition \ref{prop:rough cocylce_dyadic}, 
we obtain the following result.

\medskip

\begin{cor}
Assume that \mbox{$\big(| \cdot |_{E \otimes E}, \gamma\big)$} is exact, and let $X \colon \R \to E$ be a Wiener process based on $(E,\mathcal{H},\gamma)$, with $p$-rough path lift $\overline{\mathbf{X}} \colon \R \to T_1^{2}(E)$, for some $p \in (2,3)$. Then there exists a metric dynamical system $\big(\Omega,\mathcal{F},\P, \theta\big)$ and a $\mathcal{C}_0^{0,p-\text{var}}\big(\R,T_1^{2}(E)\big)$-valued random variable ${\mathbf{X}}$ that has the same law as $\overline{\mathbf{X}}$ and enjoys the cocycle property \ref{prop:rough cocylce_dyadic}.
\end{cor}

\bigskip

\subsection{Existence of rough driver cocycles}
\label{SubsectionExistencedrivercocycles}

Let $\big(\Omega, \mathcal{F},\P,\theta\big)$ be a measurable metric dynamical system  and $V \colon \R \times \Omega \to \mcC^0_b(\R^m,\R^m)$ a continuous time dependent random vector field. Following L.~Arnold (and others), we call $V$ a \emph{helix} if
\begin{align*}
V_{s,s+t}(\omega) = V_{0,t}(\theta_s \omega)
\end{align*}
holds for every $\omega \in \Omega$ and every $s,t \in \R$.

\begin{prop}\label{prop:existence_rough_driver_approx}
	Let $2 \leq p < 3$ and $p-2 < \rho \leq 1$ be given, and assume that ${\bfV}(\omega) = (V(\omega),\mathbb{V}(\omega))$ is a $(p,\rho)$-rough driver for every $\omega \in \Omega$ with associated vector field $W(\omega)$. Assume that $V$ is a helix. Set $\mathbb{D}_n := \{k 2^n\, :\, k \in \Z \}$, $V^n := V^{\mathbb{D}_n}$, 
	\begin{align*}
	W^n_{s,t}(x) := \frac{1}{2} \int_s^t [V^n_{du_1}, V^n_{du_2}](x)
	\end{align*}
	where the integrals are defined as Riemann-Stieltjes integrals,
	\begin{align*}
	\mathbb{V}^n_{s,t} := W^n_{s,t} + \frac{1}{2} V^n_{s,t} V^n_{s,t}
	\end{align*}
	and ${\bfV}^n :=(V^n, \mathbb{V}^n)$. Assume that 
	\begin{align*}
	d_{p,\rho}({\bfV}^n,{\bfV}) \to 0
	\end{align*}
	in probability for $n \to \infty$. Then there is an indistinguishable version of ${\bfV}$ which enjoyes the cocycle property.
\end{prop}

\medskip

\begin{Dem}
By definition of the piecewise-linear approximation, $V^n_{s,s+t}(\omega) = V^n_{0,t}(\theta_s \omega)$ for every $s,t \in \mathbb{D}_n$. Proposition \ref{prop:smooth_rp_are_cocycles} implies that also $W^n_{s,s+t}(\omega) = W^n_{0,t}(\theta_s \omega)$ for every $s,t \in \mathbb{D}_n$.	Passing to a subsequence, we may assume that $d_{p,\rho}({\bfV}^n,{\bfV}) \to 0$ for $n \to \infty$ on a subset $\tilde{\Omega}$ of full measure. This implies that the cocycle property holds for ${\bfV}$ for every $\omega \in \tilde{\Omega}$ and every $s,t \in \bigcup_n \mathbb{D}_n$, the dense set of dyadic numbers. By continuity, it then holds for every $s,t \in \R$. We can define a indistinguishable version by setting the driver equal to $0$ outside $\tilde{\Omega}$.
\end{Dem}

\medskip

\begin{example} {\normalfont
	If $M$ is a martingale vector field, there is a natural lift of $M$ to a rough driver $\mathbf{M}$, cf. \cite[Theorem 18]{BR15}. If $M$ is a martingale helix (cf. \cite[2.3.8 Definition]{Arn98}), we can apply the former proposition to the lift $\mathbf{M}$, using the approximation result given in \cite[Theorem 22]{BR15}, to obtain a indistinguishable version which is a rough driver cocycle.
	}
\end{example}

\bigskip

\subsection{Gaussian drivers}
\label{SubsectionGaussianDrivers}

The theory of stochastic flows was built in the early eighties after the pioneering works of the Russian school on the  study of the dependence of solutions to stochastic differential equations on initial conditions \cite{BL61, GS72}. One can consider as one of the early achievements of that theory the characterization of Brownian flows as solutions of stochastic differential equations driven by infinitely many Brownian motions -- this was made clear by the works of Harris, Baxendale and Le Jan \cite{Ha81, Bax80, LeJ82}. Following this line of development, we discuss in this section Gaussian drivers which take the form
\begin{align}\label{eqn:gaussian_driver_level_1}
V_{s,t}(x,\omega) = \sum_{n = 1}^{\infty} \sigma_n(x) \big(\beta^n_t-\beta^n_s\big)(\omega)
\end{align}
where $(\sigma_n)_{n\geq 1}$ is a family of vector fields in $\R^m$ and $(\beta^n)_{n\geq 1}$ is a family of continuous, real-valued Gaussian processes. By a formal calculation,
\begin{align}\label{eqn:gaussian_driver_level_2}
\begin{split}
\int_s^t [V_{du_1}, V_{du_2}](x) &= \int_s^t V^i_{s,u}(x) \partial_i V_{du}(x) -  \int_s^t \partial_i V_{s,u}(x) V^i_{du}(x) \\
&= \sum_{n,k = 1}^{\infty} \sigma_n^i(x) \partial_i \sigma_k(x) \int_s^t \beta_{s,u}^k\, d\beta_u^n - \partial_i \sigma_n(x) \sigma^i_k(x) \int_s^t \beta_{s,u}^n\, d\beta_u^k  \\
&= \sum_{n,k = 1}^{\infty} [\sigma_n, \sigma_k](x) \int_s^t \beta_{s,u}^k\, d\beta_u^n.
\end{split}
\end{align}
The calculation is formal because we have to give a meaning to the stochastic integrals (at least for $n \neq k$) and furthermore make sure that the series converge. Sufficient conditions which even lead to a rough driver are given in the next theorem. Before we state it, we recall the notion of \emph{two-dimensional $\varrho$-variation} -- see for instance \cite{FV10-2} or \cite[Chapter 10]{FH14}). If $R \colon [s,t]^2 \to \R$ is a function and $\varrho \in [1,\infty)$, set
  \begin{align*}
   \| R \|_{\varrho ; [s,t]^2} := \sup_{\mathbb{D}, \mathbb{D}'} \left( \sum_{t_i \in \mathbb{D}, t'_j \in \mathbb{D}'} \big|R(t_{i+1},t'_{j+1}) - R(t_{i},t'_{j+1}) - R(t_{i+1},t'_{j}) + R(t_{i},t'_{j})\big|^{\varrho} \right)^{\frac{1}{\varrho}} 
  \end{align*}
where the supremum ranges over all finite subsets $\mathbb{D}$ and $\mathbb{D}'$ of $[s,t]$.

\medskip

\begin{thm}
	Let $\big(\Omega, \mathcal{F},\P)$ be a probability space. Let $(\sigma_n)_{n\geq 1}$ be a family of vector fields in $\R^m$ such that
	\begin{align*}
	\sum_{n = 1}^{\infty} \| \sigma_n \|_{\mathcal{C}^{2 + \gamma}} < \infty \quad \text{and} \quad  \sum_{n, k = 1}^{\infty} \big\| [\sigma_n, \sigma_k] \big\|_{\mathcal{C}^{1 + \gamma}} < \infty
	\end{align*}
	for some $0 < \gamma \leq 1$. Furthermore, assume that there is an $\eta > 0$ and a constant $\kappa > 0$ such that
	\begin{align*}
	\sum_{n = 1}^{\infty} \big| D^i \sigma_n(x) \big| \leq \frac{\kappa}{1 + |x|^{\eta}} \quad \text{and} \quad \sum_{n, k = 1}^{\infty} \big|D^j [\sigma_n, \sigma_k](x) \big| \leq \frac{\kappa^2}{1 + |x|^{2\eta}}
	\end{align*}
	for every $x \in \R^m$, $i = 0,1,2$ and $j = 0,1$. Let $(\beta^n)_{n\geq 1}$ be a family of real-valued, centered, identically distributed, independent and continuous Gaussian processes starting at $0$ defined on the real line. Set $R(s,t) := \E\big[\beta^1_s \beta^1_t\big]$. Assume that there is a $\varrho \in \big[1,\frac{3}{2}\big)$ and a positive constant $M$ such that
	\begin{align}\label{eqn:rho_var_finite}
	\| R \|_{\varrho ; [s,t]^2} \leq M |t - s|^{1/ \varrho}
	\end{align}
	holds for every $s \leq t$. Then the following holds true.
	\begin{itemize}
		\item[\textsf{\textbf{(i)}}] The stochastic integrals
		\begin{align*}
		\int_s^t \beta_{s,u}^k\, d\beta_u^n
		\end{align*}
		exist as limit of Riemann sums in $L^q(\mathbb{P})$ for every $n,k \in \N$, $n \neq k$, $s \leq t$ and $q \geq 1$. Furthermore, the series \eqref{eqn:gaussian_driver_level_1} and \eqref{eqn:gaussian_driver_level_2} converge in $L^q(\mathbb{P})$ for any $q \geq 1$ towards continuous vector fields $V$ and $\int [V_{du_1}, V_{du_2}]$ in $\R^m$. Setting
		\begin{align*}
		W_{s,t}(x) := \frac{1}{2} \int_s^t [V_{du_1}, V_{du_2}](x) \quad \text{and} \quad \mathbb{V}_{s,t} := W_{s,t} + \frac{1}{2} V_{s,t} V_{s,t},
		\end{align*}
		the pair ${\bfV} :=(V, \mathbb{V})$ is almost surely a geometric $(p,\rho)$-rough driver for any $(p,\rho)$ satisfying $p \in (2\varrho, 3)$, $\rho \in (0, \gamma)$ and $p - 2 < \rho < 1$. Furthermore, the random variables
		\begin{align*}
			\sup_{-L \leq s \leq t \leq L} \frac{ \| V_{s,t} \|_{\mathcal{C}^{2 + \rho}}}{|t - s|^{1/p}} \qquad \text{and} \qquad \sqrt{ \sup_{-L \leq s \leq t \leq L} \frac{ \| W_{s,t} \|_{\mathcal{C}^{1+ \rho}}}{|t - s|^{2/p}}}
		\end{align*}
		have Gaussian tails for any $L \geq 0$.   \vspace{0.15cm}
		
		\item[\textsf{\textbf{(ii)}}] Assume in addition that $\big(\Omega, \mathcal{F},\P,\theta\big)$ is a measurable metric dynamical system and that all $\beta^n$ satisfy the helix property, i.e. that
		\begin{align*}
		\beta_{s,s+t}^n(\omega) = \beta^n_{0,t}(\theta_s \omega)
		\end{align*}
		holds for every $\omega \in \Omega$, $s,t \in \R$ and $n \in \N$. Then there is an indistinguishable version of  ${\bfV}$ which enjoys the cocycle property and therefore defines a rough driver cocycle.
\end{itemize}
\end{thm}

\medskip

\begin{Dem}
	We start with \textsf{\textbf{(i)}}. The condition \eqref{eqn:rho_var_finite} is the key assumption in Gaussian rough paths theory (cf. \cite{FV10-2}) and implies in particular the existence of the stochastic integrals, cf. \cite[Proposition 10.3]{FH14}. Note that it also (trivially) implies that $\| \beta^n_{s,t} \|_{L^2} \leq \sqrt{M} |t-s|^{1/(2\varrho)}$ holds for every $s \leq t$ and $n \in \N$, and Kolmogorovs continuity theorem shows that the sample paths of $\beta^n$ are $1/p$-H\"older continuous for every $p > 2\varrho$. Furthermore, the same Proposition (cf. also the discussion on p. 135 in \cite{FH14}) shows that there is a constant $C_1 > 0$ such that
	\begin{align}\label{eqn:L2_est_itint}
	\left\| \int_s^t \beta_{s,u}^k\, d\beta_u^n \right\|_{L^2(\mathbb{P})} \leq C_1 |t - s|^{1/\varrho}
	\end{align}
	holds for any $s \leq t$ and $n, k \in \N$. Fix some $L > 0$. We first claim that the series
	\begin{align}\label{eqn:gaussian_driver_level1}
	\sum_{n = 1}^{\infty} \sigma_n(x) \beta^n_t
	\end{align}
	converges uniformly in the space $\mathcal{C}\big([-L,L],\mathcal{C}_b^0(\R^m,\R^m)\big)$ almost surely. Set 
	\begin{align*}
	S_N := \sum_{n = 1}^N \| \sigma_n \|_{\mathcal{C}^0} \|\beta^n \|_{\infty;[-L,L]}
	\end{align*}
	and note that $(S_N)$ is a submartingale. Clearly,
	\begin{align*}
	\sup_{N \geq 1} \E[S_N] \leq \E\Big[ \|\beta^1 \|_{\infty;[-L,L]} \Big] \sum_{n = 1}^{\infty} \| \sigma_n \|_{\mathcal{C}^0} < \infty
	\end{align*}
	which implies, by the martingale convergence theorem, that $(S_N)$ converges almost surely. By the Weierstrass criterion for uniform convergence, the series \eqref{eqn:gaussian_driver_level1} indeed converges uniformly almost surely, and we set 
	\begin{align*}
	V_t(x) := \sum_{n = 1}^{\infty} \sigma_n(x) \beta^n_t.
	\end{align*}
	We can now repeat the same argument for the series
	\begin{align*}
	\sum_{n = 1}^{\infty} \partial_i \sigma_n(x) \beta^n_t \quad \text{and} \quad \sum_{n = 1}^{\infty} \partial^2_{i,j} \sigma_n(x) \beta^n_t
	\end{align*}
	which shows their uniform convergence and the identities
	\begin{align*}
	\partial_i V_t(x) = \sum_{n = 1}^{\infty} \partial_i \sigma_n(x) \beta^n_t, \qquad \partial^2_{i,j} V_t(x) = \sum_{n = 1}^{\infty} \partial^2_{i,j} \sigma_n(x) \beta^n_t.
	\end{align*}
	Now take $s \leq t$. We have the estimate
	\begin{align*}
	\| V_t - V_s \|_{\mathcal{C}^0} \leq | t-s|^{1/p} \sum_{n = 1}^{\infty} \| \sigma_n \|_{\mathcal{C}^0} \|\beta^n \|_{1/p-\text{H\"ol};[-L,L]}.
	\end{align*}
	Arguing as before, we can conclude that the series is almost surely convergent. The same holds for the first and second derivative, and we can conclude that
	\begin{align*}
	\sup_{-L \leq s \leq t \leq L} \frac{ \| V_{s,t} \|_{\mathcal{C}^2}}{|t - s|^{1/p}} < \infty
	\end{align*}
	almost surely. Next, our assumptions imply that there is a constant $C_2 > 0$ such that for any $q \geq 1$,
	\begin{align*}
	\big\|D^2 V_{s,t}(x) - D^2 V_{s,t}(y)\big\|_{L^q} &\leq C_2 \sqrt{q} |t-s|^{1/(2 \varrho)} |x - y|^{\gamma} \sqrt{M} \sum_{n = 1}^{\infty} \| D^2 \sigma_n \|_{\mathcal{C}^{\gamma}}, \\
	\big\| D^2 V_{s,t}(x) \big\|_{L^q} &\leq C_2 \sqrt{q} \frac{\sqrt{M} \kappa |t - s|^{1/(2\varrho)}}{1 + |x|^{\eta}},
	\end{align*}
	for any $s,t \in [-L,L]$ and $x \in \R^m$. We can now apply the Kolmogorov criterion for rough drivers \cite[Theorem 12]{BR15} to conclude that
	\begin{align*}
	\sup_{-L \leq s \leq t \leq L} \frac{ \| V_{s,t} \|_{\mathcal{C}^{2 + \rho}}}{|t - s|^{1/p}} < \infty
	\end{align*}
	almost surely for any $\rho \in (0, \gamma)$. Next we consider the series
	\begin{align*}
	\sum_{n,k = 1}^{\infty} [\sigma_n, \sigma_k](x) \int_s^t \beta_{s,u}^k\, d\beta_u^n.
	\end{align*}
	Fix $s,t \in [-L,L]$. It is easy to check that the series converges uniformly in the space $\mathcal{C}_b^0(\R^m,\R^m)$ in $L^2(\P)$, and we denote the limit by $W_{s,t}$. Using that the iterated integrals are elements in the second chaos and the $L^2$ estimates \eqref{eqn:L2_est_itint}, we can furthermore deduce that there is a constant $C_3 > 0$ such that
	\begin{align*}
	\big\|W_{s,t}(x) - W_{s,t}(y)\big\|_{L^q} &\leq C_3 q |t-s|^{1/\varrho} |x - y|^{\gamma} \sum_{n, k = 1}^{\infty} \big\| [\sigma_n, \sigma_k] \big\|_{\mathcal{C}^{1}},
	\end{align*}
	and
	\begin{align*}
	\big\|W_{s,t}(x) \big\|_{L^q} &\leq \frac{C_3 q \kappa^2 |t - s|^{1/\varrho}}{1 + |x|^{2\eta}}
	\end{align*}
	hold for every $x,y \in \R^m$, $s,t \in [-L,L]$ and $q \geq 1$. Similar estimates can be given for $V$, and we can use them for the Kolmogorov theorem \cite[Theorem 12]{BR15} to deduce that
	\begin{align*}
	\sup_{-L \leq s \leq t \leq L} \frac{ \| W_{s,t} \|_{\mathcal{C}^{0}}}{|t - s|^{2/p}} < \infty
	\end{align*}
	almost surely. For every $i = 1, \ldots, m$, it is easy to see that
	\begin{align*}
	Z_{s,t}^i = \sum_{n,k = 1}^{\infty} \partial_i[\sigma_n, \sigma_k] \int_s^t \beta_{s,u}^k\, d\beta_u^n
	\end{align*}
	exists as $L^2$-limit in the space $\mathcal{C}_b^0(\R^m,\R^m)$, and repeating the same arguments as above show that
	\begin{align*}
	\sup_{-L \leq s \leq t \leq L} \frac{ \|Z_{s,t}^i \|_{\mathcal{C}^{\rho}}}{|t - s|^{2/p}} < \infty
	\end{align*}
	almost surely for every $\rho \in (0,\gamma)$. We can now find a subsequence such that
	\begin{align*}
	\sum_{n,k = 1}^{N} [\sigma_n, \sigma_k] \int_s^t \beta_{s,u}^k\, d\beta_u^n \to V_{s,t} \qquad \text{and} \qquad
	\sum_{n,k = 1}^{N} \partial_i[\sigma_n, \sigma_k] \int_s^t \beta_{s,u}^k\, d\beta_u^n \to Z_{s,t}^i
	\end{align*}
	for $N \to \infty$ almost surely in the space $\mathcal{C}_b^0(\R^m,\R^m)$. This shows that actually
	\begin{align*}
	\sum_{n,k = 1}^{N} \partial_i[\sigma_n, \sigma_k] \int_s^t \beta_{s,u}^k\, d\beta_u^n \to \partial_i V_{s,t}
	\end{align*}
	almost surely for the same subsequence. Since we could have started with any subsequence, this shows that
	\begin{align*}
	\sum_{n,k = 1}^{\infty} \partial_i[\sigma_n, \sigma_k] \int_s^t \beta_{s,u}^k\, d\beta_u^n = \partial_i V_{s,t}
	\end{align*}  
	where the limit is understood in probability. This implies that $Z^i_{s,t} = \partial_i V_{s,t}$ almost surely for every $s,t \in [-L,L]$. We can show the same equality on set of full measure for every $s,t \in [-L,L] \cap \Q$, and continuity implies that $Z^i = \partial_i V$ almost surely. Thus we obtain
	\begin{align*}
	\sup_{-L \leq s \leq t \leq L} \frac{ \| W_{s,t} \|_{\mathcal{C}^{1+ \rho}}}{|t - s|^{2/p}} < \infty.
	\end{align*}
	Since $L > 0$ was arbitrary, this shows that $V$ and $W$ define a weakly geometric $(p,\rho)$-rough driver. They actually define a geometric $(p,\rho)$-rough driver since the piecewise-linear approximations converge, cf. part (ii). The stated tail estimates also follow from the Kolmogorov theorem \cite[Theorem 12]{BR15}.
	
\ssk	
	
For proving \textbf{\textsf{(ii)}}, note that piecewise-linear approximations of $\beta^n$ and their iterated integrals converge in $L^q(\P)$ to the original process resp. their iterated integrals, cf. \cite[Corollary 10.6 and Remark 10.7]{FH14}. We can then perform a similar strategy as above, using the Kolmogorov theorem for rough driver distance \cite[Theorem 13]{BR15}, to see that the piecewise linear approximations of $V$ and $W$ converge to the respective processes in probability. Using Proposition \ref{prop:existence_rough_driver_approx} implies the claim. We leave the details to the reader.
\end{Dem}

\medskip

\begin{example} {\normalfont
The result can be applied for $(\beta^n)$ being fractional Brownian motions with Hurst parameter $H \in (1/3,1)$ provided the $(\sigma_n)$ are sufficiently smooth and decay sufficiently fast as stated.   }
\end{example}

\bigskip

\section{Random dynamical systems induced by RDEs}
\label{SectionRDSFlows}

Let $\big(\Omega, \mathcal{F},\P,\theta\big)$ be a measurable metric dynamical system. We first consider rough differential equations on infinite dimensional spaces, studied e.g. in \cite{LCL07, Bai14, BGLY14}. Let $E$ and $F$ stand for two (non-necessarily separable) Banach spaces. Lyons' \emph{Universal Limit Theorem} states that for $\mathbf{x} \in \mathcal{C}_0^{0,p-\text{var}}\big([0,T],T_1^{\lfloor p \rfloor}(E)\big)$ and $\sigma \in L\big(E,\mathcal{C}^{\gamma}(F,F)\big)$, with $\gamma  > p$, the equation
\begin{align*}
 dy_t &=  \sigma(y_t)\, d\mathbf{x}_t, \quad y_0\in F
\end{align*}
has a unique solution path $y : [0,T]\rightarrow W$; it depends continuously on the driving rough path $\mathbf{x}$ and the initial condition $y_0$ -- see for instance \cite[Theorem 5.3]{LCL07}. Now let $\mathbf{X} \colon \R \times \Omega \to T_1^{\lfloor p \rfloor}(E)$ be a $p$-rough cocycle in the sense of Definition \ref{def:p_rough_helix}. The next statement shows the existence of a cocycle solution flow map 
$$
\varphi \colon [0,\infty) \times \Omega \times F \to F,
$$ 
for which $t \mapsto \varphi_t\big(\omega,z_0\big)$ is a solution to the equation
\begin{align}\label{eqn:rde_without_drift}
 dz_t &=  \sigma(z_t)\, d\mathbf{X}_t(\omega), \quad z_0\in F
\end{align}
on any time interval $[0,T]$.

\medskip

\begin{thm}\label{thm:RDE_inf_dim_ind_RDS}
Let $\big(\Omega, \mathcal{F},\P,\theta\big)$ be a measurable metric dynamical system and let 
$$
\mathbf{X} \colon \R \times \Omega \to T_1^{\lfloor p \rfloor}(E)
$$ 
be a $p$-rough cocycle for some $p \geq 1$. If $\sigma \in L\big(E,\mathcal{C}^{\gamma}(F,F)\big)$, for some $\gamma > p$, then there exists a unique continuous random dynamical system $\varphi$ over $\big(\Omega, \mathcal{F},\P,\theta\big)$ which solves the rough differential equation \eqref{eqn:rde_without_drift}.
\end{thm}

\medskip

\begin{Dem}
This is almost immediate. Denote by $\psi \colon [0,\infty) \times [0,\infty) \times \Omega \times F \to F$ the continuous random flow induced by equation \eqref{eqn:rde_without_drift}, obtained as a continuous function of $\bfX$. Using the characteriation of $z$ as solution to a Taylor-type numerical scheme which only involves increments of $\mathbf{X}$ \cite{Bai14, BGLY14}, one sees that $z \vert_{[s,t]}$ depends only on $z_s$ and on $(\mathbf{X}_{u,v})_{s \leq u < v \leq t}$. We therefore have
\begin{align*}
\psi\big(s + h,t + h,\omega,z_0\big) = \psi\big(s,t,\theta_h \omega, z_0\big)
\end{align*}
for every $0 \leq s \leq t$, every $h \geq 0$, $z_0 \in F$ and $\omega \in \Omega$. Using the flow property of $\psi$, it is easily seen that 
$$
\varphi\big(t,\omega,z_0\big) := \psi\big(0,t,\omega,z_0\big)
$$ 
has the cocycle property.
\end{Dem}

\medskip

The very same argument justifies that a similar result holds for rough flows associated with rough driver cocycles -- see \cite{BR15}. We record that fact here. It provides a clean extension to results on helix semimartingales reported in Arnold's book \cite{Arn98}; at the same time, its application in the setting of Gaussian rough drivers provides an interesting new class of cocycles over some non-Brownian noise space.

\medskip

\begin{thm}
Let $\big(\Omega, \mathcal{F},\P,\theta\big)$ be a measurable metric dynamical system and let $\bfV$ be a $(p,\rho)$-rough driver cocycle over $\big(\Omega, \mathcal{F},\P,\theta\big)$, for some $1\leq p<3$ and $\rho > p/3$. There exists a unique continuous random dynamical system $\varphi$ over $\big(\Omega, \mathcal{F},\P,\theta\big)$ which induces the solution flow of the equation $d\psi = {\bfV}(\psi;\,dt)$.
\end{thm}

\medskip

Let us come back to the rough differential equation \eqref{eqn:rde_without_drift}. We want to include now a drift in the equation for which we only assume a one-sided growth condition. From now on, we will work on finite dimensional spaces, that is we assume $E = \R^d$ and $F = \R^m$. We consider rough differential equations of the form
\begin{align}\label{eqn:rde_with_drift}
 dx_t &= b(x_t)\, dt + \sigma(x_t)\, d\mathbf{X}_t(\omega); \quad t \in [0,\infty)
\end{align}
where $\mathbf{X} \colon \R \times \Omega \to T_1^{\lfloor p \rfloor}(\R^d)$ is a $p$-rough cocycle, $b$ is a vector field on $\R^m$ and $\sigma \in L\big(\R^d,\mathcal{C}^{\gamma}(\R^m,\R^m)\big)$. In the following, we aim to establish existence of a random dynamical system $\varphi \colon [0,\infty) \times \Omega \times \R^m \to \R^m$ which solves the equation \eqref{eqn:rde_with_drift} in the sense of Friz-Victoir -- \cite[Section 10.3]{FV10} or \cite[Definition 1.1]{RS16}. The conditions on the drift vector field which we formulate here first appeared, in the context of rough paths, in \cite{RS16} -- see also \cite{SS16} where such conditions were used in the context of Kunita flows.   \vspace{0.15cm}

\begin{enumerate}
     \item[\textcolor{gray}{$\bullet$}] There exists a constant $C_1$ such that
     \begin{align}\label{eqn:bound_radial_growth}
      \big\langle b(x), x \big\rangle \leq C_1\big(1 + |x|^2\big), \quad \text{for every } x\in \R^m.
     \end{align}   \vspace{0.1cm}
     
     \item[\textcolor{gray}{$\bullet$}] There exists a constant $C_2$ such that
     \begin{align}\label{eqn:bound_tang_growth}
      \left| \,b(x) - \frac{\langle b(x),x \rangle \,x}{|x|^2} \right| \leq C_2\big(1 + |x|\big), \quad \text{for every } x \in \R^m \setminus \{0\}.
     \end{align}   \vspace{0.1cm}
     
     \item[\textcolor{gray}{$\bullet$}] For every $R > 0$, there exists an $R$-dependent constant $C_3(R)$ such that
     \begin{align}\label{eqn:bound_radial_growth_uniq}
     \big\langle b(x) - b(y), x-y \big\rangle \leq C_3(R)\,|x - y|^2, \quad \text{for every } x,y \in B(0,R).
     \end{align}   \vspace{0.1cm}
     
     \item[\textcolor{gray}{$\bullet$}] For every $R > 0$, there exists an $R$-dependent constant $C_4(R)$ such that
     \begin{align}\label{eqn:bound_tang_growth_uniq}
      \left|\, b(s) - b(y) - \frac{\big\langle b(x) - b(y) , x-y \big\rangle (x-y)}{| x-y |^2} \right| \leq C_4(R)\,|x-y|,
     \end{align}
     for every $x,y \in B(0,R)$, with $x \neq y$.   \vspace{0.15cm}
\end{enumerate}

It is shown in \cite[Theorem 4.3]{RS16} that equation \eqref{eqn:rde_with_drift} generates a semiflow under the assumption that the diffusivity $\sigma \in L\big(\R^d,\mathcal{C}^{\gamma + 1}(\R^m,\R^m)\big)$, for some $\gamma > p$, and that the drift $b$ is continuous and satisfies the above growth conditions. It turns out that the cocycle property of $\mathbf{X}$ also implies the existence of a cocycle $\varphi$.

\medskip

\begin{thm}
Let $\big(\Omega, \mathcal{F},\P,\theta\big)$ be a measurable metric dynamical system and let 
$$
\mathbf{X} \colon \R \times \Omega \to T_1^{\lfloor p \rfloor}(\R^d)
$$ 
be a $p$-rough cocycle, for some $p \geq 1$. Assume that $\sigma$ is a $\mathcal{C}^{\gamma + 1}(\R^m,\R^m)$-valued one form on $\R^d$, for some $\gamma > p$, and that the drift $b$ is continuous and satisfies the growth conditions \eqref{eqn:bound_radial_growth}, \eqref{eqn:bound_tang_growth}, \eqref{eqn:bound_radial_growth_uniq} and \eqref{eqn:bound_tang_growth_uniq}. Then there exists a unique continuous random dynamical system $\varphi$ over $\big(\Omega, \mathcal{F},\P,\theta\big)$ which solves the rough differential equation with drift \eqref{eqn:rde_with_drift}.
\end{thm}

\medskip

\begin{Dem}
We use the flow transformation techniques as in \cite[Theorem 4.3]{RS16}. Consider first the random flow 
$$
\psi \colon [0,\infty) \times [0,\infty) \times \Omega \times \R^m \to \R^m
$$ 
induced by the random rough differential equation
\begin{align}\label{eqn:random_rde_without_drift}
dz_t &= \sigma(z_t)\, d\mathbf{X}_t(\omega); \quad t \in [0,\infty).
\end{align}
Note that under the above regularity assumption on $\sigma$, the flow $\psi$ is differentiable \cite[Proposition 11.11]{FV10}. Set
\begin{align*}
J\big(s,t,\omega,z_0\big) := \big(D_{\xi} \psi(s,t,\omega,z_0)\big)^{-1}
\end{align*}
for $s,t \in [0,\infty)$, $z_0 \in \R^m$ and $\omega \in \Omega$. We have already seen in the proof of Theorem \ref{thm:RDE_inf_dim_ind_RDS} that we have
\begin{align*}
\psi\big(s + h,t + h,\omega,z_0\big) = \psi\big(s,t,\theta_h \omega, z_0\big)
\end{align*}
for every $0 \leq s \leq t$, $h \geq 0$, and every $\xi \in \R^m$ and $\omega \in \Omega$; the same identity holds for $J$. Let 
$$
\phi \colon [0,\infty) \times [0,\infty) \times \Omega \times \R^m \to \R^m
$$ 
denote the continuous random semiflow induced by \eqref{eqn:rde_with_drift}, whose existence and continuity is shown in \cite[Theorem 4.3]{RS16}. We aim to show that for every $0 \leq s \leq t$, every $h \geq 0$, every $x_0 \in \R^m$ and every $\omega \in \Omega$, we have
\begin{align}\label{eqn:almost_cocycle_prop_flow}
\phi\big(s+h,t+h,\omega,x_0\big) = \phi\big(s,t,\theta_h \omega,x_0\big).
\end{align}
Fix $h \geq 0$ and $0 \leq s \leq t$. Choose $T$ large enough such that $t + h \leq T$, and consider equation \eqref{eqn:rde_with_drift} on the compact time interval $[0,T]$. Fix $\omega \in \Omega$ and $x_0 \in \R^m$. In the proof of \cite[Theorem 4.3]{RS16}, it is shown that there is a $\delta > 0$ such that, for $s \leq t$, with 
\begin{align}\label{eqn:small_time_intervals}
\big\| \mathbf{X}(\theta_h \omega) \big\|_{p-\text{var};[s,t]}^p + |t - s| \leq \delta,
\end{align}
then $\phi(s,t,\theta_h \omega,\xi)$ is given by the formula
\begin{align*}
\phi\big(s,t,\theta_h \omega,x_0\big) = \psi\big(s,t,\theta_h \omega,\chi_s(t,\theta_h \omega,x_0)\big)
\end{align*}
where $u \mapsto \chi_s\big(u,\theta_h \omega,x_0\big)$ is the forward-in-time solution of the ordinary differential equation
\begin{align*}
\dot{y}_u &= J\big(s,u, \omega, y_u\big)\, b\big(\psi(s,u,\omega,y_u)\big)   \\
y_s &= x_0.
\end{align*}
Since $\big\| \mathbf{X}(\theta_h \omega) \big\|_{p-\text{var};[s,t]}^p = \big\| \mathbf{X}(\omega) \big\|_{p-\text{var};[s+h,t+h]}^p$, we also have
\begin{align*}
\phi\big(s+h,t+h, \omega,x_0\big) = \psi\big(s+h,t+h, \omega,\chi_{s+h}(t+h, \omega,x_0)\big)
\end{align*}
for the same choice of $s$ and $t$. Assume first that $0 \leq s \leq t \leq T$ are close enough to satisfy the constraint \eqref{eqn:small_time_intervals}. We claim that in this case, we have 
$$
\chi_{s + h}\big(t + h,\omega,x_0\big) = \chi_{s }\big(t , \theta_h \omega,x_0\big).
$$ 
We prove this by showing that the map $u \mapsto \chi_{s + h}\big(u + h,\omega,x_0\big)$ solves the equation
\begin{align*}
\dot{y}_u &= J\big(s,u, \theta_h\omega, y_u\big) b\big(\psi(s,u,\omega,y_u)\big)   \\
y_s &= x_0.
\end{align*}
Clearly, $\chi_{s + h}\big(s + h,\omega,x_0\big) = x_0$, and
\begin{align*}
&\frac{d}{d u} \chi_{s + h}(u + h,\omega,x_0)   \\
&= J\Big(s + h,u + h,\omega, \chi_{s + h}\big(u + h,\omega,x_0\big)\Big)\, b\Big(\psi\big(s + h,u + h,\omega,\chi_{s + h}(u + h,\omega,x_0)\big)\Big)   \\
 &= J\Big(s ,u , \theta_h \omega, \chi_{s + h}\big(u + h,\omega,x_0\big)\Big)\, b\Big(\psi\big(s ,u ,\theta_h \omega,\chi_{s + h}(u + h,\omega,x_0)\big)\Big),
\end{align*}
which implies the claim. It follows that
\begin{align*}
\phi\big(s+h,t+h,\omega,x_0\big) &= \psi\big(s + h,t + h,\omega,\chi_{s+h}(t+h,\omega,x_0)\big)   \\
&= \psi\big(s ,t ,\theta_h\omega,\chi_{s}(t, \theta_h\omega,x_0)\big) = \phi\big(s,t,\theta_h \omega,x_0\big)
\end{align*}
for $s,t$ close enough to satisfy the constraint \eqref{eqn:small_time_intervals}. Let now $s,t$ be arbitrary. Then there are numbers $(\tau_n)_{n = 0,\ldots,N}$, with $s = \tau_0 < \tau_1 < \ldots < \tau_{N-1} < \tau_N = t$,and such that every pair $(\tau_n,\tau_{n+1})$ satisfies \eqref{eqn:small_time_intervals}. We then have by the semiflow property of $\phi$
\begin{align*}
\phi\big(s,t,\theta_h \omega,x_0\big) &= \phi\big(\tau_{N-1},\tau_N,\theta_h \omega, \cdot\big) \circ \cdots \circ \phi\big(\tau_{1},\tau_0,\theta_h \omega, x_0\big)   \\
  &= \phi\big(\tau_{N-1} + h,\tau_N + h, \omega, \cdot\big) \circ \cdots \circ \phi\big(\tau_{1} + h,\tau_0 + h, \omega, x_0\big)   \\
  &= \phi\big(s + h,t + h,\omega,x_0\big).
\end{align*}
Thus we have shown that \eqref{eqn:almost_cocycle_prop_flow} indeed holds for every choice of $s \leq t$, $h$, $x_0$ and $\omega$. It is then easy to check that $\varphi\big(t,\omega,x_0\big) := \phi\big(0,t,\omega,x_0\big)$ has the cocycle property. 
\end{Dem}

\bigskip
\bigskip

\appendix

\section{Paths of bounded variation}
\label{AppendixFacts}
As usual, $(E, | \cdot |)$ denotes a real Banach space.

\medskip

\begin{definition*} {\normalfont
Let $[-a,b]$ be an interval containing $0$. The space $\mathcal{C}_0^{0,1-\text{var}}\big([-a,b], E\big)$ is defined as the closure of the set of arbitrarily often differentiable paths $x$ from $[-a,b]$ to $E$ with $x_0 = 0$ with respect to the $1$-variation norm.   }
\end{definition*}
 
\medskip

\begin{lem}\label{lemma:piecewise_diff_paths_in_C01}
The set $\mathcal{C}_0^{0,1-\text{var}}\big([-a,b], E\big)$ contains all piecewise continuously differentiable paths.
\end{lem}

\medskip

\begin{Dem}
For simplicity we assume $[-a,b] = [0,1]$. Let  $x \colon [0,1] \to E$ be piecewise differentiable, i.e. we assume that there is a dissection 
$$
\mathbb{D} = \{0 < t_1 < \ldots < t_m = 1 \}
$$ 
of $[0,1]$ such that $[t_i, t_{i+1}] \ni t \mapsto x_t $ is continuously differentiable for all $t_i$. Let $\varphi \colon \R \to [0, \infty)$ be an arbitrarily often differentiable function with support in $[-1,1]$ such that $\| \varphi \|_{L^1} = 1$. For $\varepsilon > 0$, set 
\begin{align*}
  \varphi^{(\varepsilon)}(t) := \varepsilon^{-1} \varphi\left(\frac{t}{\epsilon}\right).
\end{align*}
Using the Bochner integral, we define
\begin{align*}
  x^{(\varepsilon)}_t := \int_{\R} \varphi^{(\varepsilon)}(t - s) x_s\, ds.
\end{align*}
Lebesgue's dominated convergence theorem for the Bochner integral implies that we can differentiate under the integral, thus $x^{(\varepsilon)}$ is arbitrarily often differentiable for every positive $\varepsilon$. We claim that $\big\|x - x^{(\varepsilon)} \big\|_{1-\text{var}} \to 0$ for $\varepsilon$ converging to $0$. It is sufficient to prove that $\big\|x - x^{(\varepsilon)} \big\|_{1-\text{var};[t_i,t_{i+1}]} \to 0$ for $\varepsilon \to 0$ for all $t_i$. We have
\begin{align*}
\big\|x - x^{(\varepsilon)} \big\|_{1-\text{var};[t_i,t_{i+1}]} &= \int_{t_i}^{t_{i+1}} \big|x'_s -\big (x^{(\varepsilon)}_s\big)'\big|\, ds \leq \int_{t_i}^{t_{i+1}} \int_{\R} \big|\varphi^{(\varepsilon)}(s - u)(x'_s - x'_u) \big| \, du\, ds   \\
  &\leq \sup_{s \in [t_i,t_{i+1}]} \sup_{s - \varepsilon \leq u \leq s + \varepsilon} |x'_s - x'_u|.
\end{align*}
By uniform continuity of $x'$, the right hand side converges to $0$ for $\varepsilon \to 0$ which shows the claim.
\end{Dem}

\medskip

\begin{prop}\label{prop:C01_Polish}
The Banach space $\mathcal{C}_0^{0,1-\text{var}}\big([-a,b], E\big)$ is separable if and only if $E$ is separable.
\end{prop}

\medskip

\begin{Dem}
For simplicity, we assume that $[-a,b] = [0,1]$. Assume first that $E$ is not separable. The subset
\begin{align*} 
\big\{t \mapsto v t\, :\, v \in E \big\} \subset \mathcal{C}_0^{0,1-\text{var}}\big([0,1], E\big)
\end{align*}
is isomorphic to $E$, therefore it is not separable which implies that the space $\mathcal{C}_0^{0,1-\text{var}}\big([0,1], E\big)$ cannot be separable either. Now let $E$ be separable. Let 
$$
\mathbb{D} = \{0 < t_1 < \ldots < t_m = 1 \}
$$ 
be a dissection of $[0,1]$ and set $|\mathbb{D}| = \max_{i} |t_{i+1} - t_i|$. For a path $x \colon [0,1] \to E$, define $x^\mathbb{D} \colon [0,1] \to E$ by
\begin{align*}
  x^\mathbb{D}_t := x_{t_{i+1}} - \frac{t_{i+1} - t}{t_{i+1} - t_i} (x_{t_{i+1}} - x_{t_i}) \quad \text{for } t \in [t_i,t_{i+1}].
\end{align*}
We claim that for smooth $x \colon [0,1] \to E$, one has
\begin{align*}
  \big\| x^\mathbb{D} - x \big\|_{1-\text{var}} \to 0
\end{align*}
for $|\mathbb{D}| \to 0$. Indeed: one has
\begin{align*}
   \big\| x^\mathbb{D} - x \big\|_{1-\text{var}} = \sum_i \int_{t_i}^{t_{i+1}} \left| x'_s - \frac{x_{t_{i+1}} - x_{t_i}}{t_{i+1} - t_i} \right|\, ds
\end{align*}
and
\begin{align*}
\int_{t_i}^{t_{i+1}} \left| x'_s - \frac{x_{t_{i+1}} - x_{t_i}}{t_{i+1} - t_i} \right|\, ds \leq |t_{i+1} - t_i| \left( \sup_{t \in [t_i,t_{i+1}]} |x'_t - x'_{t_i}| + \left| x'_{t_i} - \frac{x_{t_{i+1}} - x_{t_i}}{t_{i+1} - t_i} \right| \right).
\end{align*}
By assumption, $x' \colon [0,1] \to E$ is continuous, hence uniformly continuous, so
 \begin{align*}
  \max_i \sup_{t \in [t_i,t_{i+1}]} |x'_t - x'_{t_i}| \to 0
 \end{align*}
for $|\mathbb{D}| \to 0$. For the second term, we can use the fundamental theorem of calculus for the Bochner integral \cite[Proposition A.2.3]{PR07} twice to see that
 \begin{align*}
  x_{t_{i+1}} - x_{t_i} = (t_{i+1} - t_i)\,x'_{t_i} + \int_{t_i}^{t_{i+1}} \int_{t_i}^s x''_u\, du\, ds
 \end{align*}
which implies the bound
 \begin{align*}
  \left| x'_{t_i} - \frac{x_{t_{i+1}} - x_{t_i}}{t_{i+1} - t_i} \right| \leq |t_{i+1} - t_i| \| x'' \|_{\infty}.
 \end{align*}
This shows that
 \begin{align*}
   \big\| x^\mathbb{D} - x \big\|_{1-\text{var}} \leq \max_i \sup_{t \in [t_i,t_{i+1}]} |x'_t - x'_{t_i}| + |\mathbb{D}| \| x'' \|_{\infty}
 \end{align*}
which implies the claim. Since smooth paths are dense in $\mathcal{C}_0^{0,1-\text{var}}\big([0,1], E\big)$, one can check that also $x^\mathbb{D}$ converges to $x$ in $1$-variation topology as the mesh $|\mathbb{D}|$ of the partition $\mathbb{D}$ tends to $0$, for an arbitrary $x \in \mathcal{C}_0^{0,1-\text{var}}\big([0,1], E\big)$. Now let $S \subset E$ be dense and countable and define
 \begin{align*}
  \mathcal{E} &= \big\{ x \colon [0,T] \to E\, :\, x = x^\mathbb{D} \text{ for some } \mathbb{D} = \{0 < t_1 < \ldots < t_m = 1\} \text{ with } t_i \in \Q  \\
  &\qquad \text{and } x_{t_i} \in S \text{ for all } t_i \in \mathbb{D} \big\}.
 \end{align*}
As we have the inclusion $\mathcal{E} \subset \mathcal{C}_0^{0,1-\text{var}}\big([0,1], E\big)$ from Lemma \ref{lemma:piecewise_diff_paths_in_C01}, our claim implies that $\mathcal{E}$ is dense in $\mathcal{C}_0^{0,1-\text{var}}\big([0,1], E\big)$, which shows the assertion.
\end{Dem}

\bigskip

\section{Basics on rough paths}
\label{AppendixRP}

We provide in this section some basics on rough paths; the reader is refered to the lecture notes \cite{LCL07, FH14}, or \cite{BaiLN}, for pedagogical accounts of the theory.

\medskip

We define rough path spaces for Banach space valued paths. Let $I \subset \R$ be a compact interval and $E$ a real Banach space where all tensor products are equipped with compatible norms. For a path $\mathbf{x} \colon I \to T_1^N(E)$, we define its increments setting
\begin{align*}
   \mathbf{x}_{s,t} := \mathbf{x}_{s}^{-1} \otimes \mathbf{x}_{t}
\end{align*}
for $s,t \in I$. Let $p \geq 1$ and $k \in \{1,\ldots, N\}$. The \emph{homogeneous $p$-variation metric} is defined as follows: for any two paths $\mathbf{x}, \mathbf{y} \colon I \to T_1^N(E)$, 
\begin{align*}
   d_{p-\text{var};I}(\mathbf{x},\mathbf{y}) := \max_{k = 1,\ldots,N} \left( \sup_{\mathbb{D}} \sum_{t_i \in \mathbb{D}} |\pi_k(\mathbf{x}_{t_i,t_{i+1}} - \mathbf{y}_{t_i,t_{i+1}})|^\frac{p}{k} \right)^{\frac{1}{p}}
\end{align*}
where the supremum is taken over all finite dissections $\mathbb{D}$ of the interval $I$. We also define a ``norm'' by
  \begin{align*}
   \| \mathbf{x} \|_{p-\text{var};I} := d_{p-\text{var};I}(\mathbf{x},\mathbf{1}).
\end{align*}

   Let $I \subset \R$ be a compact interval containing $0$. Note that any continuous path of bounded variation $x \colon I \to E$ has a natural lift to a path $S_N(x) \colon I \to T_1^N(E)$ given by
\begin{equation}
\label{EqDefnLift}
    \pi_k(S_N(x)_t) = \begin{cases}
                       \int_{\Delta_{[0,t]}^k} dx \otimes \cdots \otimes dx &\text{if } t \geq 0,\ 1 \leq k\leq N   \\
                       \int_{\Delta_{[t,0]}^k} dx \otimes \cdots \otimes dx &\text{if } t \leq 0,\ 1 \leq k\leq N   \\
                       1 &\text{if } k = 0.
                      \end{cases}
\end{equation}
  
\medskip   
   
\begin{definition*} {\normalfont
The space $\mathcal{C}_0^{0,p-\text{var}}\big(I,T^N_1(E)\big)$ is defined as the set of continuous paths $ \mathbf{x} \colon I \to T^N_1(E)$ with $\mathbf{x}_0 = \mathbf{1}$ for which there exists a sequence of paths $x_n \in \mathcal{C}^{0,1-\text{var}}_0(I,E)$ such that
\begin{align*}
   d_{p-\text{var};I}\big(\mathbf{x}, S_N(x_n)\big) \to 0
\end{align*}
for $n \to \infty$. We equip it with the topology induced by the metric $d_{p-\text{var};I}$. For $N = \lfloor p \rfloor$, this is the space of \emph{geometric $p$-rough paths}.   }
\end{definition*}

\medskip  
  
Note that for every $p \in [1,\infty)$ and $N \geq 1$, the natural lifting map 
$$
S_N \colon \mathcal{C}^{0,1-\text{var}}_0(I,E) \to \mathcal{C}_0^{0,p-\text{var}}\big(I,T^N_1(E)\big)
$$ 
is continuous by standard Riemann-Stieltjes estimates. If  $x \in \mathcal{C}^{0,1-\text{var}}_0(I,E)$, we call $S_{\lfloor p \rfloor}(x)$ the \emph{canonical lift} of $x$ to a rough path.

\medskip

\begin{prop}\label{prop:rp_spaces_comp_Polish}
Let $p \in [1,\infty)$ and $N \geq 1$. The space $\mathcal{C}_0^{0,p-\text{var}}\big(I,T^N_1(E)\big)$ is separable, and therefore Polish, if and only if $E$ is separable.
\end{prop}

\medskip 
 
\begin{Dem}
Assume that $\mathcal{C}_0^{0,p-\text{var}}\big(I,T^N_1(E)\big)$ is not separable. Since $S_N\big(\mathcal{C}^{0,1-\text{var}}_0(I,E)\big)$ is dense in $\mathcal{C}_0^{0,p-\text{var}}\big(I,T^N_1(E)\big)$, this set cannot be separable either. Since $S_N$ is continuous on $\mathcal{C}^{0,1-\text{var}}_0(I,E)$, this space is also not separable. Proposition \ref{prop:C01_Polish} then shows that $E$ cannot be separable. Now let $V$ be separable. From Proposition \ref{prop:C01_Polish}, there is a countable, dense subset $\mathcal{F} \subset \mathcal{C}^{0,1-\text{var}}_0(I,E)$. From continuity of $S_N$, the set $S_N(\mathcal{F}) \subset \mathcal{C}_0^{0,p-\text{var}}\big(I,T_1^N(E)\big)$ is dense, which already implies the claim.
\end{Dem}

\medskip

\section{Basics on rough drivers and rough flows}
\label{AppendixRoughFlows}

We provide in this section some elementary picture of rough drivers and rough flows, and refer the reader to the article \cite{BR15} for more information on this subject. 

\medskip

The easiest way to get a hand on these objects is probably to start from an elementary controlled ordinary differential equation. Let then $v_1,\dots,,v_\ell$ be some given smooth globally Lipschitz vector fields on $\R^m$, and $h^1,\dots,h^\ell$ be some real-valued controls, defined on some time interval $I$. The solution flow $(\varphi_{ts})_{0\leq s\leq t\leq T}$ of the controlled ordinary differential equation
$$
\dot z_t = \dot h^i_t \, v_i(z_t)
$$
enjoys then the Taylor expansion property
\begin{equation}
\label{EqTaylorFormula}
f\circ \varphi_{ts} = f + \big(h^i_t-h^i_s\big)V_if  + \left(\int_s^t\int_s^{u_1} dh^j_{u_2}\,dh^k_{u_1}\right)\,V_jV_kf + O\big(|t-s|^{>2}\big)
\end{equation}
for all smooth functions $f$. The notion of rough driver captures the essence of the different terms that appear in this local description of the dynamics.

\ssk

Let then $I$ be a compact interval. For a time-dependent vector field $V \colon I \to \mathcal{C}(\R^m,\R^m)$, set as usual $V_{s,t} := V_t - V_s$, and recall the classical one-to-one correspondence between vector fields and first order differential operators. 

\medskip

\begin{definition*} {\normalfont
Let $2\leq p<3$ and $p-2<\rho\leq 1$ be given. A \emph{weak geometric $(p,\rho)$-rough driver} is a family $\big({\bfV}_{s,t}\big)_{ s\leq t \in I}$, with 
$$
{\bfV}_{s,t} := \big(V_{s,t},\VV_{s,t}\big),
$$ 
and $\VV_{s,t}$ a second order differential operator, such that \vspace{0.1cm}
\begin{itemize}
   \item[{\sf (i)}] the vector fields $V_{s,t}$ are elements in $\mcC^{2+\rho}_b(\R^m,\R^m)$, with
$$
\underset{ s<t \in I }{\sup}\;\frac{\big\|V_{s,t}\big\|_{\mcC^{2+\rho}}}{|t-s|^{\frac{1}{p}}} < \infty,
$$ \vspace{0.1cm}

   \item[{\sf (ii)}] the second order differential operators 
   $$
   W_{s,t} := \VV_{s,t} - \frac{1}{2}V_{s,t}V_{s,t}, 
   $$
   are actually vector fields, and
   $$
   \underset{ s<t \in I}{\sup}\;\frac{\big\|W_{s,t}\big\|_{\mcC^{1+\rho}}}{|t-s|^{\frac{2}{p}}} < \infty,   
   $$ \vspace{0.1cm}

   \item[{\sf (iii)}] we have 
   $$
   \VV_{s,t} = \VV_{u,t} + V_{s,u}V_{u,t} + \VV_{s,u}, 
   $$
   for any $ s\leq u\leq t \in I$.   
\end{itemize}   }
\end{definition*}

\medskip

We use rough drivers to give a local description of the dynamics of a flow $\varphi$, under the form of a Taylor expansion formula
\begin{equation}
\label{EqTaylor}
f\circ\varphi_{ts} \simeq f + V_{ts}f + \VV_{ts}f.
\end{equation}
In the Taylor formula \eqref{EqTaylorFormula} for the controlled ordinary differential equation, the term $\big(h^i_t-h^i_s\big)V_i$ plays the role of $V_{ts}$, while the term $\left(\int_s^t\int_s^r dh^j_u\,dh^k_r\right)\,V_jV_k$ has the role of $\VV_{ts}$; check that properties {\sf (i)-(iii)} hold indeed for these two terms. We sometimes use the formal notation
\begin{align*}
W_{s,t}(x) = \int_s^t [V_{du_1}, V_{d u_2} ](x) = \frac{1}{2} \left( \int_s^t V^i_{s,u}(x) \partial_i V_{du}(x) - \int_s^t \partial_i V_{s,u}(x) V^i_{du}(x) \right).
\end{align*}
Note that if $V$ is smooth in time, we can use Riemann-Stieltjes \mbox{(or Young-)} integrals to make this notation rigorous. In particular, one can check that if $V \colon I \to \mathcal{C}^{2 + \rho}_b(\R^m,\R^m)$ is continuous with bounded variation, the above formulas for $W$ and $\VV$ define a weak geometric rough driver, called the \textit{canonical lift of} $V$.

\medskip

\begin{definition*} {\normalfont
For a weak geometric $(p,\rho)$-rough driver $\bfV$ defined on some interval $I$, set
\begin{equation*}
\label{EqDefnNorm}
\|{\bfV}\|_{p,\rho;I} := \underset{0\leq s< t\leq T}{\sup}\;\left\{\frac{\big\|V_{ts}\big\|_{\mcC^{2+\rho}}}{|t-s|^{\frac{1}{p}}} \vee \sqrt{ \;\frac{\big\|W_{ts}\big\|_{\mcC^{1+\rho}}}{|t-s|^{\frac{2}{p}}} }\right\}
\end{equation*}
and define an associated (pseudo-)metric setting
$$
d_{p,\rho;I}(\bfV,\bfV') := \big\|\bfV-\bfV'\big\|.
$$
A \emph{geometric $(p,\rho)$-rough driver} $\bfV$ is a $(p,\rho)$-rough driver for which there exists a sequence of time dependent vector fields $V_n \colon I \to \mathcal{C}^{2 + \rho}_b(\R^m,\R^m)$ which are continuous and of bounded variation in time such that their associated weak geometric rough drivers $\bfV_n$ satisfy
\begin{align*}
  d_{p,\rho;I}(\bfV,\bfV_n) \to 0
\end{align*}
for $n \to \infty$.   }
\end{definition*}

\medskip

The following definition of a solution to an equation driver by a rough driver emphasizes the fundamental role of the local description of a flow given by the Taylor-loke formula \eqref{EqTaylor}.

\medskip

\begin{definition*} {\normalfont
Given a $(p,\rho)$-rough driver $\bfV$ defined on a time interval $I$, a \emph{flow} $\varphi$ is said to \emph{solve the rough differential equation}
\begin{equation*}
d\varphi = {\bfV}(\varphi\,; dt)
\end{equation*}
if there exists a possibly $\bfV$-dependent positive constant $\delta$ such that the inequality
\begin{equation*}
\Big\|f\circ\varphi_{ts} - \Big\{f + V_{ts}f + \VV_{ts}f \Big\}\Big\|_\infty \lesssim  \|f\|_{\mcC^{2+r}}|t-s|^\frac{3}{p}
\end{equation*}
holds for all $f\in\mcC^{2+r}_b$, and all $s,t\in I$ with $0\leq t-s\leq\delta$. Such flows are called \emph{rough flows}.   }
\end{definition*}

\medskip

It can be shown that for $\rho > \frac{p}{3}$, any weak geometric $(p,\rho)$-rough driver $\bfV$ induces a flow $\psi$ of $\mathcal{C}^{\rho}$ homeomorphisms on that given time interval $I$ by solving the Kunita-type equation
\begin{align*}
 d \psi = \bfV(\psi\, ;\, dt);
\end{align*}
see \cite[Definition 4 and Theorem 5]{BR15}.

\bigskip
\bigskip

\subsection*{Acknowledgements}

I.B. thanks the U.B.O. for their hospitality. SR and MS gratefully acknowledge the support by the DFG Research Unit FOR 2402.

\bigskip
\bigskip

\bibliographystyle{alpha}
\bibliography{refs}

\bigskip
\bigskip

\noindent \textcolor{gray}{$\bullet$} {\sf I. Bailleul} - {\small Institut de Recherche Mathematiques de Rennes, France.}   \vspace{0.1cm}

\vspace{-0.1cm}\noindent \hspace{2.2cm}{\it ismael.bailleul@univ-rennes1.fr}   \vspace{0.3cm}

\noindent \textcolor{gray}{$\bullet$} {\sf S. Riedel} - {\small Institut f\"ur Mathematik, Technische Universit\"at Berlin, Germany.}   \\
\vspace{-0.1cm}\noindent \hspace{2.2cm}{\it riedel@math.tu-berlin.de}   \vspace{0.3cm}

\noindent \textcolor{gray}{$\bullet$}  {\sf M. Scheutzow} - {\small Institut f\"ur Mathematik, Technische Universit\"at Berlin, Germany.}   \\ 
\vspace{-0.1cm}\noindent \hspace{2.2cm}{\it ms@math.tu-berlin.de}

\end{document}